\documentclass[a4paper,11pt]{amsart}

\usepackage{graphicx}
\usepackage{amsmath}
\usepackage{amssymb}
\usepackage{amsthm}
\usepackage{amscd}
\usepackage{fullpage}
\usepackage[all]{xy}
\usepackage{tikz-cd}
\usepackage{latexsym}
\usepackage{mathrsfs}
\usepackage{dsfont} 
\usepackage{hyperref}

\newtheorem{theorem}{Theorem}[section]
\newtheorem{lemma}[theorem]{Lemma}
\newtheorem{corollary}[theorem]{Corollary}
\newtheorem{proposition}[theorem]{Proposition}
\newtheorem{introtheorem}{Theorem}

\theoremstyle{definition}
\newtheorem{remark}[theorem]{Remark}

\newtheorem{definition}[theorem]{Definition}
\newtheorem{observation}[theorem]{Observation}
\newtheorem{notation}[theorem]{Notation}
\newtheorem{setup}[theorem]{Setup}

\newcommand{\A}{\mathcal{A}}
\newcommand{\C}{\mathcal{C}}
\newcommand{\D}{\mathcal{D}}
\newcommand{\E}{\mathcal{E}}
\newcommand{\I}{\mathcal{I}}
\newcommand{\K}{\mathbb{K}}
\newcommand{\R}{\mathbb{R}}
\renewcommand{\S}{\Sigma}
\newcommand{\T}{\mathcal{T}}
\newcommand{\X}{\mathcal{X}}
\newcommand{\Z}{\mathbb{Z}}
\newcommand{\Zc}{\mathcal{Z}}
\renewcommand{\a}{\alpha}
\renewcommand{\b}{\beta}
\newcommand{\g}{\gamma}
\newcommand{\io}{\iota}
\newcommand{\ka}{\kappa}
\renewcommand{\l}{\lambda}
\newcommand{\om}{\omega}
\renewcommand{\t}{\tau}
\newcommand{\emp}{\varnothing}
\newcommand{\lora}{\longrightarrow}
\newcommand{\ovl}{\overline}
\newcommand{\ovs}{\overset}
\newcommand{\unl}{\underline}

\DeclareMathOperator{\Hom}{Hom}
\DeclareMathOperator{\End}{End}
\DeclareMathOperator{\ind}{ind}
\DeclareMathOperator{\Rep}{Rep}
\DeclareMathOperator{\rep}{rep}
\DeclareMathOperator{\Ker}{Ker}
\DeclareMathOperator{\Coker}{Coker}
\DeclareMathOperator{\Vect}{Vect}
\DeclareMathOperator{\Tr}{Tr}
\DeclareMathOperator{\thick}{thick}
\DeclareMathOperator{\pwf}{pwf}
\DeclareMathOperator{\Proj}{Proj}
\DeclareMathOperator{\Inj}{Inj}
\DeclareMathOperator{\rad}{rad}

\hypersetup{colorlinks,
linkcolor = {red!50!black},
citecolor = {blue!50!black},
urlcolor = {blue!80!black}
}

\newcommand{\harxiv}[1]{ \href{http://arxiv.org/abs/#1}{\texttt{arXiv:#1}}}

\begin{document}
\title{Negative Calabi--Yau discrete cluster categories\\ via Nakayama representations and persistence theory}

\author{Sofia Franchini} 
\email{s.franchini1@lancaster.ac.uk}

\begin{abstract}
We introduce infinite discrete versions of the symmetric Nakayama representations by using techniques of persistence theory. After stabilising, we obtain a family triangulated categories which can be regarded as negative Calabi--Yau versions of the Igusa--Todorov discrete cluster categories of type $A$. We describe their geometric model and AR theory.
\end{abstract}

\maketitle
{\small
\tableofcontents}

\addtocontents{toc}{\protect\setcounter{tocdepth}{1}}

\section{Introduction}

Calabi--Yau triangulated categories (CY for short) satisfy an important duality on their $\Hom$-spaces: for each pair of objects there is an isomorphism $\Hom(a,b)\cong D\Hom(b,\S^w a)$ which is natural in $a$ and $b$, where $\S$ is the shift functor and $w\in\Z$ stands for the CY parameter. In representation theory CY triangulated categories are related, for instance, to cluster theory, quivers with potentials, and symmetric algebras. The classical cluster categories are $2$-CY; from a quiver with potential we can obtain a $3$-CY triangulated category; the perfect derived category of a symmetric algebra is $0$-CY, while its stable module category is $(-1)$-CY.

CY triangulated categories are difficult to study in general, therefore it is useful to find classes of examples to develop intuition. The Igusa--Todorov discrete cluster categories and the Holm--J{\o}rgensen categories provide tractable and interesting classes of CY triangulated categories well suited for this purpose. For a positive integer $m$, Igusa and Todorov defined a cluster category $\C_m$ having a geometric model consisting of an $\infty$-gon $\Zc_m$ with $m$ two-sided accumulation points in \cite{IT}. The category $\C_m$ is an infinite discrete version of the classical cluster category $\C(A_n)$ of type $A_n$ introduced in \cite{BMRRT}. Igusa--Todorov categories can also be regarded as generalisations of the Holm--J{\o}rgensen category $\T_2$, first introduced and studied in \cite{HJ1, J}; indeed $\C_1$ and $\T_2$ are triangle equivalent. The categories $\T_2$ and $\C_m$ (and their completions) are the subject of intensive study, see for instance \cite{ACFGS, CKP,CZZ, CG, F, Fr, GHJ, GZ, HJ1, LP,  M1,  M2, N, PY}. Both $\C_m$ and $\T_2$ are $2$-CY. 

There exist different CY versions of the Holm--J{\o}rgensen category $\T_2$, which are denoted by $\T_w$, where $w\in\Z$ stands for the CY parameter. These categories were first considered in \cite{HJ2, HJY} and have nice geometric models, introduced in \cite{C1, HJ2}. Our aim for this paper is to find a $(-1)$-CY version of $\C_m$ which also is a ``higher number of accumulation points" version of $\T_{-1}$. We denote our desired category by $\C_{-1,m}$. 

Defining the category $\C_{-1,m}$ can be considered as a first step for subsequently defining also lower CY versions of $\C_m$. Negative CY triangulated categories are the natural counterpart of triangulated categories having a positive CY dimension. They are interesting because they show a dual homological behaviour: in \cite{ZZ} the authors proved that in the negative case there are only trivial t-structures, while the converse holds for the positive case. Negative CY triangulated categories also provide an environment for studying simple-minded systems, which are the negative CY analogue of cluster-tilting subcategories, see for instance \cite{CP2}. For the categories $\T_w$ with $w\leq -1$, the torsion pairs, (co-)t-structures, Riedtmann configurations, simple-minded systems, and their mutation theory were studied in \cite{C1, C2, CP1, HJY}.

It is not possible to adapt the Igusa--Todorov and Holm--J{\o}rgensen constructions in order to obtain $\C_{-1,m}$. To define $\C_{-1,m}$ we combine Nakayama representations and techniques of persistence theory. By \cite[Theorem 7.0.5]{A}, the $(-1)$-CY orbit category $\D^b(A_n)/\S^2\t$, where $\t$ is the Auslander--Reiten translation, is triangle equivalent to the stable category of representations $\rep(C_n, \rad^{n+1})$ where $C_n$ is the oriented cycle with $n$ vertices and $\rad^{n+1}$ is the set of relations generated by the paths of length $n+1$. 

The bound path algebra $\K C_n/\rad^{n+1}$ is an example of Nakayama algebra, i.e. a finite dimensional algebra for which the indecomposable left and right modules have a unique composition series. Nakayama algebras are equivalent to the categories of representations $\rep(Q,I)$ where $Q$ is a quiver of type $A_n$ with linear orientation or is the oriented cycle $C_n$. Rock and Zhu defined the categories of \emph{continuous Nakayama representations}, denoted by $\Rep^{\pwf}(\R,\ka)$ and $\Rep^{\pwf}(S^1,\ka)$ for the linear and cyclic case respectively, in \cite{RZ}. To do so, they considered certain representations of $\R$ (also called \emph{persistence modules}) and of $S^1$, subject to some conditions imposed by a \emph{(pre)-Kupisch function} $\ka$. The representations of $\R$ and and of $S^1$ are important in topological data analysis and persitence theory, and were studied in \cite{C,HR} from a representation theoretic point of view; while $\ka$ represents the continuous analogue of the Kupisch series for Nakayama algebras and establishes the length for the indecomposable projective objects.

In \cite{RZ} Rock and Zhou also proved that the category of modules over a Nakayama algebra can be regarded as an abelian subcategory of $\Rep^{\pwf}(S^1,\ka)$ for a suitable $\ka$. Inspired by their method, we introduce the category of \emph{infinite discrete symmetric Nakayama representations}, which we denote by $\rep(\Zc_m,\ka)$. We first define a certain Kupisch function $\ka$, and we then consider those representations of $\Rep^{\pwf}(S^1,\ka)$ which are ``constant" in between any two consecutive elements of $\Zc_m$ and in the neighbourhoods of the accumulation points (see Definition \ref{definition infinite discrete symmetric nakayama representations}). We obtain the following result.

\begin{introtheorem}[{Theorems \ref{theorem rep is wide} and \ref{theorem decomposition theorem for rep}, Proposition \ref{proposition hom spaces}}]
The category $\rep(\Zc_m,\ka)$ is a Krull--Schmidt abelian subcategory of $\Rep^{\pwf}(S^1,\ka)$. The isoclasses of its indecomposable objects are in bijection with certain intervals of $\R$. The $\Hom$-spaces between indecomposable objects are at most two-dimensional and can be understood from the intersections of the corresponding intervals.
\end{introtheorem}

We recall that $\rep(C_n,\rad^{n+1})$ is a Frobenius category and $P_i = I_i$ for each vertex $i$ of $C_n$. We obtain that $\rep(\Zc_m,\ka)$ has these same properties. We also prove that, while each indecomposable object of $\rep(C_n,\rad^{n+1})$ has a unique (finite) composition series, each indecomposable object of $\rep(\Zc_m,\ka)$ has a unique, possibly countable infinite, chain of subobjects for which the quotient of any two consecutive elements is simple. We regard this as a (possibly countable infinite) composition series. Therefore, $\rep(\Zc,\ka)$ is an infinite version of $\rep(C_n,\rad^{n+1})$.

\begin{introtheorem}[{Proposition \ref{proposition enough injective objects} and Theorem \ref{theorem uniserial}}]	
The category $\rep(\Zc_m,\ka)$ has enough projectives and enough injectives. For each $z\in\Zc_m$ we have that $P_z = I_z$, i.e. the indecomposable projective object and the indecomposable injective object at $z$ coincide. Each indecomposable object of $\rep(\Zc_m,\ka)$ admits a unique, possibly countable infinite, composition series.
\end{introtheorem}

We define $\C_{-1,m}$ as the stable category of $\rep(\Zc_m,\ka)$ and we obtain a geometric model in terms of the $\infty$-gon $\Zc_m$. Finally, we prove that the category $\C_{-1,m}$ is both an analogue of the Igusa--Todorov category $\C_m$ and a generalisation of the Holm--J{\o}rgensen category $\T_{-1}$.

\begin{introtheorem}[{Proposition \ref{proposition bijection intervals and -1-admissible arcs}, Corollary \ref{corollary -1-cy}, and Theorem \ref{theorem triangle equivalence HJ category}}]
The category $\C_{-1,m}$ is a Krull--Schmidt $(-1)$-CY triangulated category. The isoclasses of its indecomposable objects are in bijection with certain arcs of $\Zc_m$. The $\Hom$-spaces between indecomposable objects are at most one-dimensional and can be understood from the relative position of the corresponding arcs. Moreover, $\C_{-1,1}$ is triangle equivalent to $\T_{-1}$.
\end{introtheorem}

\subsection*{Acknowledgments}
The author thanks David Pauksztello for useful mathematical discussions, for the suggestions, and for carefully reading previous drafts of this paper. The author thanks Eric Hanson for useful mathematical discussion. The author acknowledges financial support from the EPSRC through the grant EP/V050524/1.

\section{Background}

For the rest of this paper, we fix a field $\K$. Let $\C$ be an additive category. We say that $\C$ is \emph{$\K$-linear} if its $\Hom$-sets are $\K$-vector spaces, and $\C$ is \emph{$\Hom$-finite} if its $\Hom$-spaces are finite dimensional. Whenever we consider a subcategory $\X$ of $\C$, we assume it to be full and we sometimes write $\X\subseteq\C$. We say that $\X$ is an \emph{additive subcategory} of $\C$ if $\X$ contains the zero object, is closed under isomorphisms, finite direct sums, and direct summands. We say that $\C$ is \emph{Krull--Schmidt} if each object decomposes into a finite direct sum of objects each having local endomorphism ring. By \cite[Theorem 4.2]{K}, this decomposition is unique, up to isomorphism and reordering the summands. We denote by $\ind\C$ the class of indecomposable objects of $\C$. 

\subsection{Frobenius categories}\label{section frobenius categories}

Let $\A$ be an abelian category, and $\E$ be an additive subcategory of $\A$.

\begin{itemize}
\item We say that $\E$ is \emph{extension-closed}, or \emph{closed under extensions}, if for each short exact sequence $0\lora e_1\lora a\lora e_2\lora 0$ with $e_1,e_2\in\E$ we have that $a\in\E$. 
\item We say that $\E$ is a \emph{wide subcategory} of $\A$ if $\E$ is closed under kernels, cokernels, and extensions. Note that in particular $\E$ is an abelian subcategory of $\A$.
\item Assume that $\E$ is closed under extensions. We say that the category $\E$ together with the short exact sequences of $\A$ which restrict to $\E$ is an \emph{exact subcategory} of $\A$. We denote by $\Proj\E$ and $\Inj\E$ the classes of projective and injective objects of $\E$ respectively. We refer for instance to \cite[Defintion 2.1]{B} for a definition of exact category which does not rely on an ambient abelian category. 
\item We say that $\E$ is a \emph{Frobenius subcategory} of $\A$ if $\E$ is an exact subcategory of $\A$, $\E$ has enough projectives and enough injectives, i.e. each object of $\E$ has a projective cover and an injective envelope in $\E$, and $\Proj\E = \Inj\E$.
\end{itemize}

From a Frobenius subcategory $\E$ of $\A$, we obtain the \emph{stable category} $\unl\E$. This category has a triangulated structure induced by the exact structure of $\E$. We recall the definition of $\unl\E$ below, for more details we refer to \cite[Section 5.1.2]{Z} for instance. These facts will be useful in Section \ref{section a negative calabi-yau triangulated category} for establishing the properties of the category $\C_{-1,m}$.

The category $\unl\E$ has as objects the same objects of $\E$. To define the morphisms of $\unl\E$, we first introduce an equivalence relation $\sim$ on the $\Hom$-sets of $\E$. Given $a,b\in\E$, we denote the class of morphisms of $\E$ which factor through a projective object of $\E$ by $\Proj(a,b)$. Given morphisms of $\E$ $f,g\colon a\to b$, we have that $f\sim g$ if and only if $f-g\in\Proj(a,b)$. We denote by $[f]$ the equivalence class of $f$ and we say that $[f]$ is the \emph{stabilisation} of $f$. The morphisms of $\unl\E$ are exactly of the form $[f]$.

We now discuss the triangulated structure on $\unl\E$. Given $a\in\E$, consider its injective envelope $\io\colon a\to i$ and extend it to a short exact sequence $0\lora a\ovs{\io}\lora i\lora \S a\lora 0$  where $\S a  = \Coker\io$. This assignment induces an automorphism $\S\colon\unl\E\to\unl\E$ which is the shift functor of $\unl\E$. Consider an exact sequence $0\lora a\ovs{f}\lora e\ovs{g}\lora b\lora 0$, since $i$ is injective, there exists $\phi\colon e\to i$ such that $\phi f = \io$. Moreover, by the universal property of the cokernels, there exists $h\colon b\to \S a$ such that the following diagram commutes.
\[
\begin{tikzcd}
0\ar[r] & a\ar[r,"f"]\ar[d,"1"] & e\ar[r,"g"]\ar[d,"\phi"] & b\ar[r]\ar[d,"h"] & 0\\
0\ar[r] & a\ar[r,"\io"] & i\ar[r] & \S a\ar[r] & 0
\end{tikzcd}
\]
The isoclasses of the sequences of the form $a\ovs{[f]}\lora e\ovs{[g]}\lora b\ovs{[h]}\lora \S a$ give the class of triangles of $\unl\E$. 

Note that if $\E$ is $\K$-linear, $\Hom$-finite, and Krull--Schmidt, then so is $\unl\E$.

\subsection{Irreducible morphisms and almost split triangles}

Throughout this section we assume $\E$ to be a $\K$-linear, $\Hom$-finite, Krull--Schmidt Frobenius category, and $\T = \unl\E$. 

Let $f\colon a\to b$ be a morphism in an additive category. We recall that $f$ is a \emph{split monomorphism} if there exists $f'\colon b\to a$ such that $f'f = 1_a$. The morphism $f$ is \emph{left almost split} if $f$ is not a split monomorphism and each $g\colon a\to c$ which is not a split monomorphism factors through $f$. The notions of \emph{split epimorphism} and of \emph{right almost split} morphism are dual. The morphism $f$ is \emph{irreducible} if $f$ is not a split monomorphism and not a split empimorphism, and if $f = hg$ then $g$ is a split monomorphism or $h$ is a split epimorphism. 

Each irreducible morphism of $\T$ is obtained by stabilising an irreducible morphism of $\E$, and each almost split triangle in $\T$ is obtained by stabilising an almost split sequence in $\E$. We refer to \cite{ASiSk} and \cite{H} for some background on AR theory. We believe that these facts are well known, but we could not find references in the literature, therefore we provide arguments with Proposition \ref{proposition irreducible and left right almost split} and Corollary \ref{corollary almost split triangles after stabilising} for the convenience of the reader. 

\begin{proposition}\label{proposition irreducible and left right almost split}
Let $\E$ be a $\Hom$-finite, $\K$-linear, Krull--Schmidt Frobenius abelian category and $f\colon a\to b$ be a morphism in $\E$. The following statements hold.
\begin{enumerate}
\item If $a,b\in\ind\E\setminus\Proj\E$, then $f$ is irreducible in $\E$ if and only if $[f]$ is irreducible in $\unl\E$.
\item If $f$ is a monomorphism and $a\in\ind\E\setminus\Proj\E$, then $f$ is left almost split if and only if $[f]$ is left almost split.
\item  If $f$ is an epimorphism and $b\in\ind\E\setminus\Proj\E$, then $f$ is right almost split if and only if $[f]$ is right almost split.
\end{enumerate}
\end{proposition}
\begin{proof}
We prove (1). Assume that $f\colon a\to b$ is irreducible in $\E$. Since $f$ is not a split monomorphism nor a split epimorphism, $[f]$ is not a split monomorphism nor a split epimorphism. Assume that $[f] = [h][g] = [hg]$ for some $g\colon a\to c$ and $h\colon c\to b$ in $\E$. We have that $f-hg = p_2p_1$ for some $p_1\colon a\to q$ and $p_2\colon q\to b$ where $q\in\Proj\E = \Inj\E$. Since $f$ is irreducible and $\E$ is abelian, $f$ is either a monomorphism or an epimorphism. Assume that $f$ is a monomorphism, the other case is dual. Since $q\in\Inj\E$, there exists $\a\colon b\to q$ such that $\a f = p_1$. As a consequence, $hg = (1_b-p_2\a)f$. Since $\End_{\E}(b)$ is local, $1_b-p_2\a$ is an isomorphism because $p_2\a$ is not an isomorphism. Indeed, if $p_2\a$ is an isomorphism then $\a$ is a split monomorphism and $b\in\Proj\E$, giving a contradiction. Thus, $f = (1_b-p_2\a)^{-1}hg$ and, from the fact that $f$ is irreducible, $g$ is a split monomorphism or $(1_b-p_2\a)^{-1}h$ is a split epimorphism. It is straightforward to check that this implies that $[g]$ is a split monomorphism or $[h]$ is a split epimorphism. This proves that $[f]$ is irreducible.
	
Now assume that $a,b\in\ind\E\setminus\Proj\E$ and that $[f]\colon a\to b$ is irreducible in $\unl\E$. We have that $f$ is not a split monomorphism nor a split epimorphism, otherwise $[f]$ is. If $f = hg$ for some $g\colon a\to c$ and $h\colon c\to b$, then $[f] = [h][g]$. Therefore, $[g]$ is a split monomorphism or $[h]$ is a split epimorphism, and then $g$ is a split monomorphism or $h$ is a split epimorphism. We conclude that $f$ is irreducible in $\E$.
	
Now we prove (2), the proof of (3) is dual. If $f$ is left almost split, then $f$ is not a split monomorphism and therefore $[f]$ is not a split monomorphism in $\unl\E$. Consider $[g]\colon a\to c$ in $\unl\E$ which is not a split monomorphism, and then $g\colon a\to c$ is not a split monomorphism in $\E$. Since $f$ is left almost split, there exists $h\colon b\to c$ such that $hf = g$, and as a consequence $[h][f] = [g]$. This proves that $[f]$ is left almost split.
	
Now assume that $[f]$ is left almost split. We show that $f$ is left almost split. Since $[f]$ is not a split monomorphism, $f$ is not a split monomorphism. Now consider a morphism $g\colon a\to c$ which is not a split monomorphism in $\E$. Then $[g]\colon a\to c$ is not a split monomorphism, and then there exists $[h]\colon b\to c$ such that $[h][f] = [g]$, i.e. $hf = g+p_2p_1$ for some $p_1\colon a\to q$ and $p_2\colon q\to c$ with $q\in\Proj\E = \Inj\E$. Since $f$ is a monomorphism, there exists $\a\colon b\to q$ such that $\a f = p_1$, and then $(h-p_2\a)f = g$. We conclude that $f$ is left almost split.
\end{proof}

We recall that an exact sequence $0\lora a\ovs{f}\lora e\lora b\lora 0$ in $\E$ is an \emph{almost split sequence}, or an \emph{Auslander--Reiten (AR) sequence}, if $a,b\in\ind\E$ and $f$ is left almost split. The object $a$ is denoted by $\t(b)$. The following corollary of Proposition \ref{proposition irreducible and left right almost split} relates the almost split triangles of $\unl\E$ to the almost split sequences of $\E$.

\begin{corollary}\label{corollary almost split triangles after stabilising}
Let $\E$ be a $\Hom$-finite, $\K$-linear, Krull--Schmidt Frobenius category. If $0\lora a\ovs{f}\lora e\ovs{g}\lora b\lora 0$ is an almost split sequence in $\E$, then $a\ovs{[f]}\lora e\ovs{[g]}\lora b\lora\S a$ is an almost split triangle in $\unl\E$.
\end{corollary}

\subsection{Continuous representations}

In this section we discuss the representations over $\R$ and the circle $S^1$, and the continuous Nakayama representations. These will be useful in Section \ref{section infinite discrete symmetric nakayama representations} for defining the infinite discrete symmetric Nakayama representations. For more details on the material presented here we refer to \cite{C, HR, RZ}.

\subsubsection{Representations of $\R$}\label{section representations of R} 

We regard set of real numbers $\R$ as a category: the objects of $\R$ are the real numbers, and for any $s,t\in\R$, if $s\leq t$ there is a unique morphism $f_{st}\colon s\to t$. For each $t\in\R$ the morphism $f_{tt}$ coincides with the identity $1_t$. 

We denote by $\Vect\K$ the category of vector spaces over $\K$. A \emph{representation} of $\R$ over $\K$ is a covariant functor $M\colon \R\to \Vect\K$. A \emph{morphism} of representations is a natural transformation. We denote by $\Rep\R$ the category of representations of $\R$. A representation $M\in\Rep\R$ is \emph{pointwise finite} if $\dim M(t)<\infty$ for each $t\in\R$. We denote by $\Rep^{\pwf}\R$ the category of pointwise finite representations of $\R$.

We fix some notation and terminology. Given $M,N\in\Rep\R$, by $\Hom_{\R}(M,N)$ we indicate, with an abuse of notation, the set of morphisms $M\to N$ in the category $\Rep\R$. For an interval $U\subseteq\R$, the \emph{interval representation} $M_U$ is given by
\[
M_U(t) = \begin{cases}
\K & \text{if $t\in U$,}\\
0 & \text{otherwise}
\end{cases}
\]
with $M_U(f_{st})\colon M_U(s)\to M_U(t)$ equal to $1_{\K}$ if $s,t\in U$ and $s\leq t$, and equal to $0$ otherwise.

\begin{definition}[{\cite[Definition 2.3]{RZ}}]
Let $U,V\subseteq\R$ be intervals. The \emph{left intersection} of $V$ and $U$ is
\[
V\cap_L U = 
\begin{cases}
V\cap U & \text{if $v<u$ for any $(v,u)\in ((V\setminus U)\times U)\cup(V\times (U\setminus V))$,}\\
\emp & \text{otherwise.}
\end{cases}	
\]
\end{definition}

\begin{figure}[h]
\centering
\includegraphics[height = 2cm]{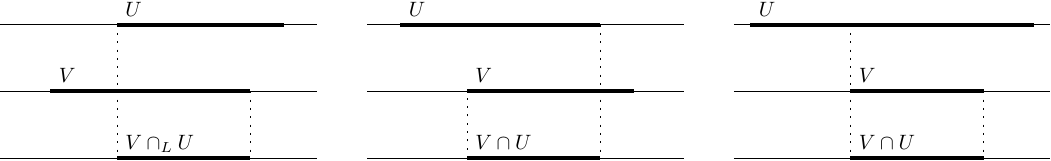}
\caption{Illustration of left intersections of the intervals $U$ and $V$. On the left, $V\cap_L U = V\cap U\neq\emp$. In the center and on the right, $V\cap U\neq\emp$ but $V\cap_L U = \emp$.}
\label{figure left intersection}
\end{figure}

We refer to Figure \ref{figure left intersection} for an illustration. It is straightforward to check that $\Hom_{\R}(M_U,M_V)\cong\K$ if $V\cap_L U\neq\emp$, and $\Hom_{\R}(M_U,M_V) = 0$ otherwise.

\begin{definition}\label{definition standard morphisms}
Let $U,V\subseteq\R$ be bounded intervals. We define the \emph{standard morphism} $\phi\colon M_U\to M_V$ as $\phi(t) = 1_{\K}$ for each $t\in V\cap_L U$ and $\phi(t) = 0$ for each $t\in\R\setminus(V\cap_L U)$.
\end{definition}

Note that the composition of two standard morphisms is still a standard morphism. 

The following theorem gives the decomposition of the representations of $\R$.

\begin{theorem}[{\cite[Theorem 1.1]{C}}] 
Any pointwise finite representation of $\R$ decomposes uniquely, up to isomorphism and reordering the summands, as a direct sum of possibly infinitely many interval representations.
\end{theorem}

\subsubsection{Representations of $S^1$}\label{section representations of circle}

Similarly as for $\R$, we regard the circle $S^1$ as a category. Given $x,y\in S^1$ with $x\neq y$, there is a path $g_{xy}\colon x\to y$ moving from $x$ to $y$ along $S^1$ in the anticlockwise direction. By concatenating $g_{xy}$ and $g_{yx}$, we obtain a path $\om_x = g_{yx}g_{xy}\colon x\to x$ which consists of moving from $x$ to $x$ in the anticlockwise direction around $S^1$ exactly once. By convention, $\om_x^0$ and $g_{xx}$ coincide with the lazy paths $x\to x$, and we often denote them by $1_x$. The objects of $S^1$, as a category, are the points of the circle, and for each pair of points $x,y\in S^1$ the morphisms $x\to y$ are the paths of the form $g_{xy}\om_x^n = \om_y^n g_{xy}$ with $n\in\Z$.

Analogously as for $\R$, a \emph{representations} of $S^1$ over $\K$ is a covariant functor $M\colon S^1\to \Vect \K$. We denote by $\Rep S^1$ the category of representations of $S^1$, and by $\Rep^{\pwf} S^1$ the category of pointwise finite representations. If $M,N\in\Rep S^1$, by $\Hom_{S^1}(M,N)$ we denote, with an abuse of notation, the set of morphisms $M\to N$ in the category $\Rep S^1$.

\begin{figure}[h]
\centering
\includegraphics[height = 5cm]{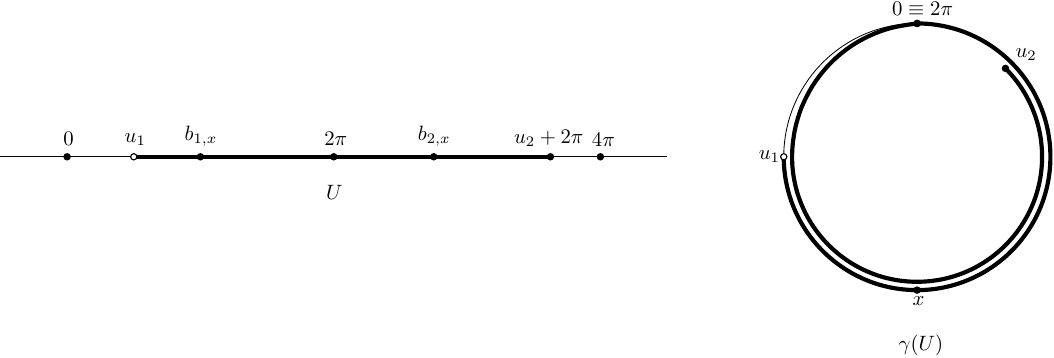}
\caption{The bounded interval $U = (u_1,u_2+2\pi]\subseteq\R$ and the cyclically ordered set $\g(U)\subseteq S^1$. The empty circles indicate that the endpoints are excluded, while the black circles indicate that the endpoints are included.}
\label{figure string representation}
\end{figure}

We introduce some notation from \cite{RZ}. Consider a bounded interval $U\subseteq\R$ and the covering map $\g\colon \R \to S^1, t \mapsto t\mod 2\pi$. Let $x\in S^1$ and consider the set $\g^{-1}(x)\cap U$, we denote its cardinality by $n_x$ and its elements by $b_{1,x}<\cdots<b_{n_x,x}$. Note that $n_x = 0$ if and only if $x\notin \g(U)$. Figure \ref{figure string representation} provides an illustration. We now introduce some representations which are building blocks of pointwise finite representations.

\begin{definition}[{\cite[Section 2.2]{RZ}}]\label{definition string}
Let $U\subseteq\R$ be a bounded interval. Keeping the notation above, the representation $\ovl{M}_U$ of $S^1$ defined as follows is called a \emph{string}. 
\begin{itemize}
\item For each $x\in S^1$ we have that $\ovl{M}_U(x) = \K b_{1,x}\oplus\dots\oplus \K b_{n_x,x}$.
\item For each $x,y\in S^1$ we have that
\[
\ovl{M}_U(g_{xy})(b_{i,x}) =  
\begin{cases}
b_{j,y} & \text{if there exists $b_{j,y}$ such that $0\leq b_{j,y}-b_{i,x}<2\pi$,}\\
0 & \text{otherwise.}
\end{cases} 
\]
\end{itemize}
\end{definition}

The representation $\ovl{M}_U$ is pointwise finite because $U$ is bounded. Moreover, $\ovl{M}_U(\om_x)(b_{i,x}) = b_{i+1,x}$ if $i\neq n_x$, and $\ovl{M}_U(\om_x)(b_{n_x,x}) = 0$.

\begin{theorem}[{\cite[Theorem 3.8]{HR}, \cite[Corollaries 2.8, 2.9]{RZ}}]\label{theorem indecomposable strings}
Let $U,V\subseteq\R$ be bounded intervals. The following statements hold.
\begin{enumerate} 	
\item We have that $\ovl{M}_U\cong \ovl{M}_V$ if and only if $U = V+2n\pi$ for some $n\in\Z$.
\item The ring $\End_{S^1}(\ovl{M}_U)$ is local, and therefore $\ovl{M}_U$ is indecomposable.
\item We have that $\Hom_{S^1}(\ovl{M}_U,\ovl{M}_V) \cong \bigoplus_{n\in\Z}\Hom_{\R}(M_U,M_{V+2n\pi})$.
\end{enumerate}
\end{theorem}

Note that, since the intervals $U$ and $V$ are bounded, the direct sum in the statement above is finite. The following theorem provides a decomposition of the pointwise finite representations of $S^1$. We refer to \cite[Definition 3.5]{HR} for the definition of \emph{band}.

\begin{theorem}[{\cite[Theorem 5.6]{HR}}]\label{theorem decomposition into strings and bands}
A pointwise finite representation of $S^1$ decomposes uniquely, up to isomorphism and reordering the summands, as a direct sum of possibly infinitely many strings and finitely many bands.
\end{theorem}

Standard morphisms in $\Rep\R$ induce morphisms between strings in $\Rep S^1$. We introduce some notation which will be useful for Definition \ref{definition morphisms of strings induced by standard morphisms}.

\begin{figure}[h]
\centering
\includegraphics[height = 2cm]{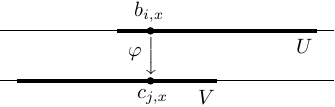}
\caption{Illustration of Notation \ref{notation definition morphisms of strings induced by standard morphisms}.}
\label{figure phi}
\end{figure}

\begin{notation}\label{notation definition morphisms of strings induced by standard morphisms}
Let $U,V\subseteq\R$ be bounded intervals, $\phi\colon M_U\to M_V$ be the standard morphism, and $x\in S^1$ be such that $\ovl M_U(x)\neq 0$. We have that $\ovl M_U(x) = \bigoplus_{i = 1}^{n_x}\K b_{i,x}$ where $\{b_{i,x}\}_{i = 1}^{n_x} = \g^{-1}(x)\cap U$. Note that if $b_{i,x}\in V\cap_L U$ for some $i$, then there exists a unique $c_{j,x}\in\g^{-1}(x)\cap V$ such that $b_{i,x} = c_{j,x}$, see Figure \ref{figure phi}. With a light abuse of notation, we denote $\phi(b_{i,x}) = c_{j,x}$ if $b_{i,x}\in V\cap_L U$, and $\phi(b_{i,x}) = 0$ otherwise.
\end{notation}

\begin{definition}[{\cite[p. 44]{RZ}}]\label{definition morphisms of strings induced by standard morphisms}
Keeping $U$, $V$, and $\phi$ as in Notation \ref{notation definition morphisms of strings induced by standard morphisms}, the morphism $\ovl\phi\colon \ovl M_U\to \ovl M_V$ is defined as $\ovl\phi(x)(b_{i,x}) = \phi(b_{i,x})$ for each $x\in S^1$ such that $\ovl M_U(x)\neq 0$.
\end{definition} 

The following lemma will be useful for proving the factorisation properties of the morphisms of the category of infinite discrete Nakayama representations, which will be defined in Section \ref{section infinite discrete symmetric nakayama representations} (see Propositions \ref{proposition factorisation properties in rep} and \ref{proposition zero compositions in rep}).

\begin{lemma}\label{lemma morphisms in rep}
Let $U,V, W\subseteq\R$ be bounded intervals. The following statements hold.
\begin{enumerate}
\item Let $\phi\colon M_U\to M_V$ be a standard morphism. Then  $\ovl\phi = 0$ if and only if $\phi = 0$, or equivalently $V\cap_L U = \emp$.
\item If $\phi\colon M_U\to M_V$ and $\psi\colon M_V\to M_W$ are standard morphisms, then $\ovl{\psi\phi} = \ovl\psi\hspace{0.05 cm} \ovl\phi$.
\item Let $\phi\colon M_U\to M_V$ and $\psi\colon M_{U+2n\pi}\to M_{V+2n\pi}$ be standard morphisms and $n\in\Z$. Then $\ovl\phi = \ovl\psi$.
\item Assume that there exists a unique $n\in\Z$ such that $(V+2n\pi)\cap_L U\neq \emp$. Then any non-zero morphism $f\colon \ovl{M}_U\to \ovl M_V$ is of the form $f = \l\ovl\phi$, where $\phi\colon M_U\to M_{V+2n\pi}$ is the standard morphism and $\l\in \K^*$ .
\end{enumerate} 
\end{lemma}
\begin{proof}
We prove statement (1). If $\phi = 0$ then $\ovl\phi = 0$ by Definition \ref{definition morphisms of strings induced by standard morphisms}. If $\phi\neq 0$ then there exists $t\in V\cap_L U\neq\emp$. Now consider $x = \g(t)\in S^1$, the sets $\g^{-1}(x)\cap U\neq \emp$ and $\g^{-1}(x)\cap V$ are non-empty because $t$ belongs to both. Since $t\in\g^{-1}(x)\cap U = \{b_{i,x}\}_{i = 1}^{n_x}$, we have that $t = b_{i,x}$ for some $i$. Therefore, $b_{i,x}\in V\cap_L U$, $\ovl\phi(x)(b_{i,x}) = \phi(b_{i,x})\neq 0$, and as a consequence $\ovl\phi\neq 0$.

Now we prove (2). If $\phi = 0$ or $\psi = 0$, then $\psi\phi = 0$ and $\ovl{\psi\phi} = 0$. Moreover, $\ovl\phi = 0$ or $\ovl\psi = 0$ and then $\ovl\psi\hspace{0.05 cm}\ovl\phi = 0 = \ovl{\psi\phi}$. Now assume that $\phi\neq 0$ and $\psi \neq 0$, i.e. $V\cap_L U\neq\emp$ and $W\cap_L V\neq\emp$. Let $x\in S^1$ be such that $\ovl M_U(x)\neq 0$, and consider $b_{i,x}\in\g^{-1}(x)\cap U$. We have the following possibilities: $W\cap_L U\neq\emp$ or $W\cap_L U = \emp$. Assume that the first case holds, the second case is similar. If $b_{i,x}\in W\cap_L U$, then $b_{i,x}\in V\cap_L U$, and the following equalities hold 
\[
\ovl\psi(x)\ovl\phi(x)(b_{i,x}) = \ovl\psi(x)(\phi(b_{i,x})) = \ovl\psi(x)(c_{j,x}) = \psi(c_{j,x}) = \psi\phi(b_{i,x}) = \ovl{\psi\phi}(x)(b_{i,x})
\] 
where $c_{j,x}\in\g^{-1}(x)\cap V$ is such that $c_{j,x} = b_{i,x}$. Now assume that $b_{i,x}\notin W\cap_L U$, then $\ovl{\psi\phi}(x)(b_{i,x}) = \psi\phi(b_{i,x}) = 0$. Moreover, if $\phi(b_{i,x}) = 0$ then $\ovl\psi(x)(\phi(b_{i,x})) = 0$, and if $\phi(b_{i,x}) \neq 0$, i.e. $b_{i,x}\in U\cap_L V$, then $\ovl\psi(x)(\phi(b_{i,x})) = \psi\phi(b_{i,x}) = 0$ because $b_{i,x}\notin W\cap_L U$. In both cases $\ovl\psi(x)\ovl\phi(x)(b_{i,x}) = \ovl\psi(x)(\phi(b_{i,x})) = 0 = \ovl{\psi\phi}(x)(b_{i,x})$. Therefore, we obtain that $\ovl\psi\hspace{0.05 cm}\ovl\phi = \ovl{\psi\phi}$.

Statement (3) is straightforward and follows from  Definition \ref{definition string}. We prove statement (4). By Theorem \ref{theorem indecomposable strings} we have that $\Hom_{S^1}(\ovl M_U,\ovl M_V)\cong\K$ and, since $(V+2n\pi)\cap_L U\neq \emp$, we have that $\phi\neq 0$ and then $\ovl\phi\neq 0$. Thus, $f = \l\ovl\phi$ for some $\l\in\K^*$.
\end{proof}

\subsubsection{Continuous Nakayama representations}

The continuous Nakayama representations of $S^1$ are pointwise finite representations of $S^1$ which satisfy some conditions determined by a map called a \emph{Kupisch function}. This map is the continuous analogue of the Kupisch series for Nakayama algebras and determines the length of the projective representations.

\begin{definition}[{\cite[Definition 3.9]{RZ}}]
A \emph{Kupisch function} $\ka\colon \R\to \R^{>0}$ is a map such that
\begin{enumerate}
\item $\ka(t+2\pi) = \ka(t)$ for each $t\in\R$, and
\item $t_1+\ka(t_1)\leq t_2+\ka(t_2)$ for each $t_1,t_2\in\R$ such that $t_1\leq t_2$.
\end{enumerate}
\end{definition} 

\begin{definition}[{\cite[Definition 3.9]{RZ}}]\label{definition continuous nakyama representations}
Let $\ka$ be a Kupisch function and $M\in\Rep^{\pwf}S^1$. We say that $M$ is a \emph{continuous Nakayama representation} of $(S^1,\ka)$ if each indecomposable direct summand of $M$ is a string of the form $\ovl{M}_U$, where $U\subseteq \R$ is a bounded interval such that $U\subseteq[\inf U,\inf U+\ka(\inf U)]$. The category of continuous Nakayama representations of $(S^1,\ka)$ is denoted by $\Rep^{\pwf}(S^1,\ka)$.
\end{definition}

\begin{remark}[{\cite[Remark 3.10]{RZ}}]\label{remark continuous nakayama representations}
The category $\Rep^{\pwf}(S^1,\ka)$ is an abelian subcategory of $\Rep(S^1)$. Any $M\in\Rep^{\pwf}(S^1,\ka)$ decomposes uniquely, up to isomorphism and reordering the summands, as the direct sum of possibly infinitely many strings.
\end{remark}

\subsection{The $\infty$-gon $\mathcal{Z}_m$}\label{section the ifinity-gon}

\begin{figure}[h]
\centering
\includegraphics[height = 4cm]{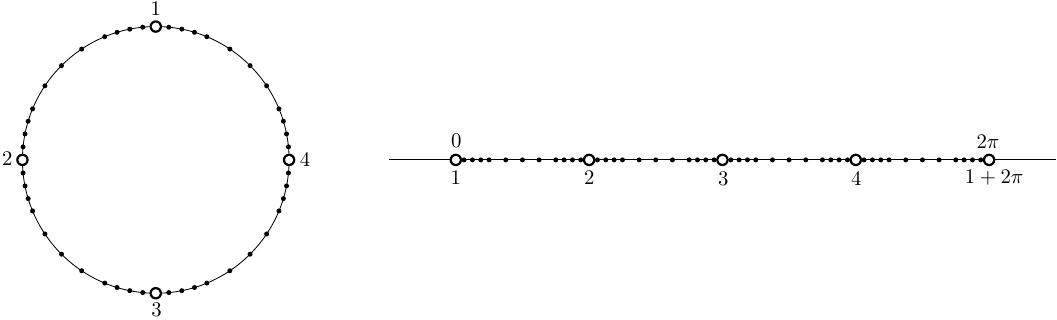}
\caption{The $\infty$-gon $\Zc_4$ regarded as a subset of $S^1$ and of the interval $(0,2\pi)\subseteq\R$.}
\label{figure infinity-gon}
\end{figure}

We now describe the geometric object used for studying the properties of discrete cluster categories. Given a positive integer $m$, we denote by $\Zc_m$ the subset of the unit circle $S^1$ obtained by embedding $m$ copies of $\Z$ in $\Zc_m$. Each copy of $\Z$ is denote by $\Z^{(p)}$ for $p\in[m]$. The elements of $\Zc_m$ accumulate to $m$ two-sided accumulation points, which do not belong to $\Zc_m$. We index the accumulation points of $\Zc_m$ by $1,\dots, m$ and we denote by $[m]$ the set of accumulation points $\{1,\dots,m\}$. We often regard $\Zc_m$ as a subset of the interval $(0,2\pi)\subseteq\R$, where the accumulation point $1\in [m]$ and the real number $0\in\R$ are identified. In particular, $(0,2\pi)$ induces a linear order on $\Zc_m$. It will be clear from the context if we regard $\Zc_m$ as a subset of $S^1$ or as a subset of $\R$. We denote the successor of $z\in\Zc_m$ by $z^+$ and the predecessor by $z^-$. We refer to Figure \ref{figure infinity-gon} for an illustration of $\Zc_m$. 

\section{Infinite discrete symmetric Nakayama representations}\label{section infinite discrete symmetric nakayama representations}

In this section we introduce the infinite discrete symmetric Nakayama representations. Given a finite dimensional Nakayama algebra $A$, in \cite[Section 6]{RZ} the authors defined a Kupisch function $\ka_A$ which encodes the length of the indecomposable projective objects of the module category $\mathrm{mod}\hspace{0.1cm} A$. They then proved that $\mathrm{mod}\hspace{0.1cm} A$ embeds into $\Rep^{\pwf}(S^1,\ka_A)$. We follow a similar approach: we define a Kupisch function $\ka$ from the $\infty$-gon $\Zc_m$, and we obtain the category of representations of $(\Zc_m,\ka)$ as an abelian subcategory of $\Rep^{\pwf}(S^1,\ka)$. The following is an infinite version of \cite[Definition 6.1]{RZ} when $A = \K C_n/\rad^{n+1}$ and $C_n$ is the oriented cycle with $n$ vertices.

Note that each $s\in[0,2\pi)$ is either an accumulation point of $\Zc_m$, i.e. $s\in[m]$, or $s\in[z,z^+)$ for some $z\in\Zc_m$.

\begin{definition}
Let $s,t\in\R$ be such that $t \equiv s\mod 2\pi$ with $s\in[0,2\pi)$. We define	
\[
\ka(t) = 
\begin{cases}
2\pi + z^+ -s & \text{if $s\in[z,z^+)$ for some for some $z\in\Zc_m$,}\\
2\pi & \text{if $s\in[m]$.}	
\end{cases}
\]
\end{definition}

Note that the map $\ka\colon \R\to \R^{>0}$ is $2\pi$-periodic.

\begin{lemma}
The map $\ka$ is a Kupisch function.
\end{lemma}
\begin{proof}
Since $\ka$ is $2\pi$-periodic, it remains to check that for any $t_1,t_2\in\R$ if $t_1\leq t_2$ then $t_1+\ka(t_1)\leq t_2+\ka(t_2)$. Without loss of generality we can assume that $t_1\in [0, 2\pi)$ and $t_2\in [0,2\pi)+2h\pi$ for some $h\in\Z$ with $h\geq 0$.

If $h = 0$ then we can proceed with a case analysis where we distinguish the cases when $t_1$ and $t_2$ are accumulation points of $\Zc_m$ or not. We show the case where $t_1\in[m]$ and $t_2\in [z,z^+)$ for some $z\in\Zc_m$. Since $t_1< t_2<z^+$, we have that $t_1+\ka(t_1) = t_1+2\pi<z^++2\pi = t_2+\ka(t_2)$. The other cases are straightforward.

If $h = 1$, let $s = t_2-2\pi\in[0,2\pi)$ and note that either $t_1\leq s$ or $t_1>s$. If $t_1\leq s$ then, from the case $h = 0$ above, $t_1+\ka(t_1)\leq s+\ka(s)\leq s+2\pi+\ka(s) = t_2+\ka(t_2)$. If $t_1>s$, then we can divide the proof into cases where $t_1$ and $s$ belong or not to $[m]$. We show the case where $t_1\in [z,z^+)$ for some $z\in\Zc_m$ and that $s\in[m]$. Since $z^+<2\pi<s+2\pi$, we have that $t_1+\ka(t_1) = 2\pi+z^+ < s+4\pi = t_2+2\pi = t_2+\ka(t_2)$. The other cases are straightforward.

If $h = 2$, then $t_2-2\pi\in[2\pi,4\pi)$. From the case $h = 1$, we have that $t_1+\ka(t_1)\leq t_2-2\pi+\ka(t_2-2\pi)  = t_2-2\pi+\ka(t_2) <t_2+\ka(t_2)$. If $h\geq 3$ we can proceed analogously. We conclude that $t_1+\ka(t_1)\leq t_2+\ka(t_2)$. 
\end{proof}

We now define the category $\rep(\Zc_m,\ka)$.

\begin{definition}\label{definition infinite discrete symmetric nakayama representations} 
Let $M\in\Rep^{\pwf}(S^1,\ka)$. We say that $M$ is an \emph{infinite discrete symmetric Nakayama representation} if it satisfies the following conditions.
\begin{enumerate}
\item If $x,y\in S^1$ are such that $z< x\leq y \leq z^+$ for some $z\in\Zc_m$, then $M(g_{xy})\colon M(x)\to M(y)$ is an isomorphism.
\item For each $p\in[m]$ there exist $a,b\in \Zc_m$ such that $a<p<b$ are in cyclic order, and if $x,y\in S^1$ are such that $a<x\leq y\leq b$ are in cyclic order, then $M(g_{xy})\colon M(x)\to M(y)$ is an isomorphism.
\end{enumerate}
We denote by $\rep(\Zc_m,\ka)$ the full subcategory of $\Rep^{\pwf}(S^1,\ka)$ of infinite discrete symmetric Nakayama representations.
\end{definition}

Given $M\in\rep(\Zc_m,\ka)$, we observe that, for each $x\in S^1$, $M(x)$ can be obtained by extending the action of $M$ at the points of $\Zc_m$. The same holds for the connecting morphisms $g_{xy}$. We can therefore consider the objects of $\rep(\Zc_m,\ka)$ to be ``discrete". 

With Proposition \ref{proposition projective injective objects} we will prove that at each $z\in\Zc_m$ the indecomposable projective object $P_z$ coincides with the indecomposable injective object $I_z$. Moreover, in Section \ref{section simple objects rep} we will prove that each indecomposable object of $\rep(\Zc_m,\ka)$ has a unique composition series. 

We now prove that $\rep(\Zc_m,\ka)$ is a wide subcategory of $\Rep^{\pwf}(S^1,\ka)$ and is therefore abelian. First, we have the following lemma.

\begin{lemma}\label{lemma rep is wide}
Let $L,M,N\in\Rep^{\pwf}(S^1,\ka)$, $f\colon L\to M$ and $g\colon M\to N$ be morphisms, and $x,y\in S^1$. The following statements hold.
\begin{enumerate}
\item Assume that $M\cong L\oplus N$. Then $M(g_{xy})$ is an isomorphism if and only if so are $L(g_{xy})$ and $N(g_{xy})$.
\item If $0\lora L\lora M\lora N\lora 0$ is a short exact sequence in $\Rep^{\pwf}(S^1,\ka)$ and $L(g_{xy})$ and $N(g_{xy})$ are isomorphisms, then so is $M(g_{xy})$.
\item If $L(g_{xy})$ and $M(g_{xy})$ are isomorphisms, then so is $\Ker f(g_{xy})$.
\item If $M(g_{xy})$ and $N(g_{xy})$ are isomorphisms, then so is $\Coker g(g_{xy})$.
\end{enumerate}
\end{lemma}
\begin{proof}
We prove (1). Since $M\cong L\oplus N$ and $M(g_{xy})$ is an isomorphism, $L(g_{xy})\oplus N(g_{xy})$ is an isomorphism. Thus, both $L(g_{xy})$ and $N(g_{xy})$ are isomorphisms. The other implication is similar. The statements (2), (3), and (4) follow, respectively, from the Five Lemma and the universal properties of kernel and cokernel.
\end{proof}

\begin{theorem}\label{theorem rep is wide}
The category $\rep(\Zc_m,\ka)$ is a wide subcategory of $\Rep^{\pwf}(S^1,\ka)$.
\end{theorem}
\begin{proof}
We prove that $\rep(\Zc_m,\ka)$ is closed under extensions. For checking that $\rep(\Zc_m,\ka)$ is an additive subcategory of $\Rep^{\pwf}(S^1,\ka)$ closed under kernels and cokernels, we can proceed similarly. Let $0\lora L\lora M\lora N\lora 0$ be a short exact sequence in $\Rep^{\pwf}(S^1,\ka)$ with $L,N\in\rep(\Zc_m,\ka)$, we prove that $M\in\rep(\Zc_m,\ka)$. 

Let $x,y\in S^1$ be such that $z<x\leq y\leq z^+$ for some $z\in\Zc_m$. Thus, $L(g_{xy})$ and $N(g_{xy})$ are isomorphisms and, by Lemma \ref{lemma rep is wide}, $M(g_{xy})$ is an isomorphism. Now consider $p\in[m]$. There exist $a_L,b_L\in\Zc_m$ such that $a_L<p<b_L$ are in cyclic order and $L(g_{xy})$ is an isomorphism for each $x,y\in S^1$ such that $a_L\leq x\leq y\leq b$. Moreover, there exist $a_N,b_N\in\Zc_m$ with the same property for $N$. Therefore, there exist $a', b'$ such that both $L(g_{xy})$ and $N(g_{xy})$ are isomorphisms for each $a'< x\leq y\leq b'$. By Lemma \ref{lemma rep is wide} we obtain that $M(g_{xy})$ is an isomorphism for each $a'< x\leq y\leq b'$. We conclude that $M\in\rep(\Zc_m,\ka)$. 	
\end{proof}

\section{Indecomposable objects}

This section is devoted to establish the Krull--Schmidt decomposition theorem for the category $\rep(\Zc_m,\ka)$. We also describe its indecomposable objects and we arrange them into a coordinate system which gives the AR quiver of $\rep(\Zc_m,\ka)$.

\subsection{The Krull--Schmidt property}

We define the following set of intervals.
\[
\I = \{(u_1,u_2+2h\pi]\subseteq\R\mid  u_1,u_2\in\Zc_m, h\in\Z,\text{ and } u_1^+\leq u_2+2h\pi\leq u_1^++2\pi\}
\]
Note that either $h = 0$ or $h = 1$. The condition $u_1^+\leq u_2+2h\pi\leq u_1^++2\pi$ ensures that the ``shortest" intervals of $\I$ are of the form $(u_1,u_1^+]$ and that the ``longest" intervals are of the form $(u_1,u_1^++2\pi]$. In Sections \ref{section projective-injective objects rep} and \ref{section simple objects rep} we will prove that these intervals correspond, respectively, to the simple objects and to the indecomposable projective-injective objects of $\rep(\Zc_m\,\ka)$.

With the following proposition we prove that conditions (1) and (2) of Definition \ref{definition infinite discrete symmetric nakayama representations} imply that the isoclasses of indecoposable objects of $\rep(\Zc_m,\ka)$ are in bijection with the intervals in $\I$.

Note that the intervals of $\I$ are similar to the intervals of \cite[Section 6]{RZ} which were used for describing the indecomposable objects of $\mathrm{mod}\hspace{0.1cm}A$ regarded as an abelian subcategory of $\Rep(S^1,\ka_A)$, where $A$ is a Nakayama algebra.

\begin{proposition}\label{proposition indecomposable objects of rep}
Let $M\in\ind\Rep^{\pwf}(S^1,\ka)$. Then $M\in\rep(\Zc_m,\ka)$ if and only if $M\cong \ovl{M}_{U}$ for some $U\in\I$.
\end{proposition}
\begin{proof}
Assume that $M\cong\ovl M_U$ for some $U = (u_1,u_2+2h\pi]\in\I$, we check that $\ovl M_U\in\rep(\Zc_m,\ka)$. Since $u_1,u_2\in\Zc_m$, $\ovl M_U(g_{xy})$ is an isomorphism for each $x,y\in S^1$ such that $z<x\leq y\leq z^+$ are in cyclic order for some $z\in\Zc_m$. Moreover, for each $p\in[m]$ there exist $a,b\in\Zc_m$ such that $\ovl M_U(g_{xy})$ is an isomorphism for each $a<x\leq y\leq b$ in cyclic order. We refer to Figure \ref{figure proposition strings} for an illustration. We conclude that $M\cong\ovl M_U\in\rep(\Zc_m,\ka)$.	
	
\begin{figure}[h]
\centering
\includegraphics[height = 5cm]{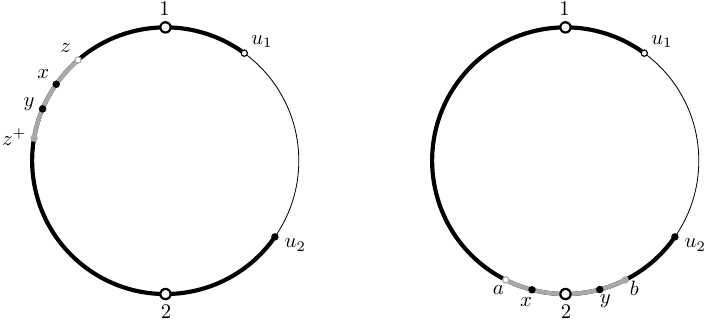}
\caption{Illustration of the argument of Proposition \ref{proposition indecomposable objects of rep} in $\Zc_2$.}
\label{figure proposition strings}
\end{figure}
	
Now assume that $M\in\rep(\Zc_m,\ka)$. By Definition \ref{definition continuous nakyama representations}, $M\cong \ovl{M}_U$ for some bounded interval $U\subseteq\R$ such that $U\subseteq [\inf U,\inf U+\ka(\inf U)]$. Without loss of generality we can assume that $\inf U\in[0,2\pi)$. Indeed, by Theorem \ref{theorem indecomposable strings}, $\ovl M_U\cong \ovl M_{U+2n\pi}$ for each $n\in\Z$. Since $0\leq\inf U\leq \sup U\leq \inf U+\ka(\inf U)\leq 4\pi$, we have that $\sup U\in[0,4\pi)$. We write $\inf U = u_1$ and $\sup U = u_2+2h\pi $ with $u_2\in [0,2\pi)$ and $h\in\{0,1\}$. If $u_1,u_2\in [m]$, then we have a contradiction with condition (2) of Definition \ref{definition infinite discrete symmetric nakayama representations} . If $u_1,u_2\notin\Zc_m$, then we have a contradiction with condition (1). Thus, $u_1,u_2\in\Zc_m$. In particular, since $\emp\neq U\subseteq[u_1,u_1+\ka(u_1)] = [u_1,u_1^++2\pi]$, we have that $u_1^+\leq u_2+2h\pi\leq u_1^++2\pi$. Finally, $u_1\notin U$ and $u_2+2h\pi\in U$, otherwise we have a contradiction with condition (1) of Definition \ref{definition infinite discrete symmetric nakayama representations}. We conclude that $U = (u_1,u_2+2h\pi]\in\I$.
\end{proof}

We recall that, from Theorem \ref{theorem decomposition into strings and bands}, each object of $\Rep(S^1,\ka)$ decomposes as a possibly infinite direct sum of indecomposable objects. With Proposition \ref{proposition finite direct summands} we prove that condition (2) of Definition \ref{definition infinite discrete symmetric nakayama representations} implies that the objects of $\rep(\Zc_m,\ka)$ have only finitely many indecomposable direct summands.

\begin{proposition}\label{proposition finite direct summands}
Let $M\in\rep(\Zc_m,\ka)$. Then $M$ has only finitely many indecomposable direct summands.
\end{proposition}
\begin{proof}
Since $\rep(\Zc_m,\ka)$ is closed under direct summands, the indecomposable direct summands of $M$ are objects of $\rep(\Zc_m,\ka)$. For each $p\in[m]$ we fix $a_p<p<b_p$ such that $M(g_{xy})$ is an isomorphism for each $x,y\in S^1$ such that $a_p<x\leq y\leq b_p$ are in cyclic order. Such  $a_p$ and $b_p$ exist by Definition \ref{definition infinite discrete symmetric nakayama representations}. Now consider the set $\Zc_m' = \Zc_m\setminus\bigcup_{p\in[m]}\{x\in S^1\mid a_p<x<b_p \text{ are in cyclic order}\}$ and note that it is finite. 
	
By Proposition \ref{proposition indecomposable objects of rep} each indecomposable direct summand of $M$ is of the form $\ovl M_U$ with $U = (u_1,u_2+2h\pi]\in\I$. We prove that each indecomposable direct summand $\ovl M_U$ of $M$ is such that $u_2\in\Zc_m'$. If $u_2\notin \Zc_m'$ then there exist $p\in[m]$ and $x,y\in S^1$ such that $a_p<u_2<b_p$ are in cyclic order. Let $x = u_2$ and $y\in S^1$ be such that $a_p<x\leq y\leq b_p$ are in cyclic order, then $\ovl M_U(g_{xy})$ is not an isomorphism, see Figure \ref{figure finite direct summands}. By Lemma \ref{lemma rep is wide} this contradicts the fact that $M(g_{xy})$ is an isomorphism. Therefore $u_2\in\Zc_m'$.
	
\begin{figure}[h]
\centering
\includegraphics[height = 5cm]{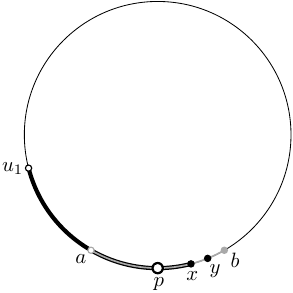}
\caption{Illustration of the argument of Proposition \ref{proposition finite direct summands}.}
\label{figure finite direct summands}
\end{figure}
	
Now consider $z\in\Zc_m'$. There are only finitely many direct summands of $M$ of the form $\ovl M_U$ with $U = (u_1,u_2+2h\pi]$ and $u_2 = x$, otherwise $\dim M(x) = \infty$. Since $\Zc_m'$ is a finite set, we conclude that $M$ has only finitely many indecomposable direct summands.
\end{proof}

Finally, we have the following decomposition theorem.

\begin{theorem}\label{theorem decomposition theorem for rep}
The category $\rep(\Zc_m,\ka)$ is Krull--Schmidt. Moreover, the indecomposable objects of $\rep(\Zc_m,\ka)$ are exactly those of the form $\ovl M_U$ with $U\in\I$, up to isomorphism.
\end{theorem}

\subsection{The coordinate system}

\begin{figure}[h]
\centering
\includegraphics[height = 5cm]{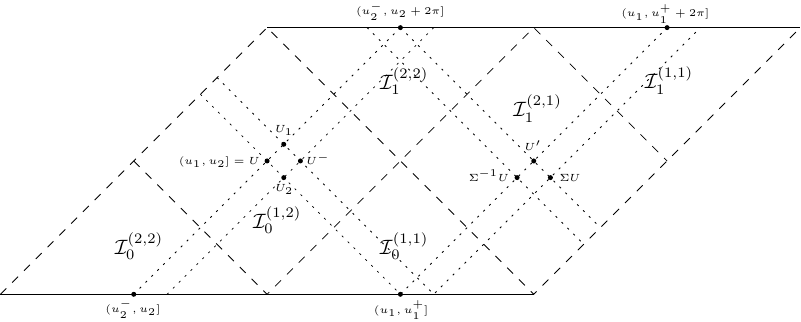}
\caption{The coordinate system of $\rep(\Zc_2,\ka)$.}
\label{figure ar quiver rep}
\end{figure}

For each $p,q\in[m]$ and $h\in\{0,1\}$ we introduce the set of intervals
\[
\I^{(p,q)}_h = \left\{(u_1,u_2+2h\pi]\in\I\mid u_1\in\Z^{(p)}\text{ and } u_2\in\Z^{(q)}\right\}.
\]
We can arrange the intervals of $\I$, or equivalently the isoclasses of indecomposable objects of $\rep(\Zc_m,\ka)$, into a coordinate system having
\begin{itemize}
\item $2m$ components of type $\Z A_{\infty}$, each corresponding to the sets $\I^{(p,p)}_h$ for $p\in[m]$ and $h\in\{0,1\}$, and
\item $2\binom{m}{2}$ components of type $\Z A_{\infty}^{\infty}$, each corresponding to the sets $\I^{(p,q)}_h$ for $p,q\in[m]$, $p\neq q$, and $h\in\{0,1\}$.
\end{itemize}
In Section \ref{section irreducible morphisms and almost split sequences in rep} we will prove that this coordinate system gives the AR quiver of $\rep(\Zc_m,\ka)$. We now introduce some intervals of $\I$. We refer to Figure \ref{figure ar quiver rep} for an illustration of the coordinate system and of the intervals defined below.

\begin{definition}\label{definition intervals}
Let $U = (u_1,u_2+2h\pi]\in\I$, we define the following intervals.	\begin{align*}
&U_1  = (u_1^-,u_2+2h\pi] & &\S^{-1}U  = (u_2,u_1^++2(1-h)\pi]\\
&U_2  = (u_1,u_2^-+2h\pi] & &\S U  = (u_2^-,u_1+2(1-h)\pi]\\
&U^-  = (u_1^-,u_2^-+2h\pi] & &U'  = (u_2^-,u_1^++2(1-h)\pi]
\end{align*}
\end{definition}

The intervals above play a central role in $\rep(\Zc_m,\ka)$. With Proposition \ref{proposition factorisation through projectives} we will prove that $U'$ determines a set of intervals $V\in\I$ for which any non-zero morphism $\ovl M_U\to \ovl M_V$ factors through a projective-injective object of $\rep(\Zc_m,\ka)$. The intervals $U_1$, $U_2$, and $U^-$ determine the almost split sequence starting at $\ovl M_U$, see Proposition \ref{proposition almost split sequences in rep}. After stabilising $\rep(\Zc_m,\ka)$, $\S U$ and $\S^{-1} U$ determine the actions of the shift functor and its inverse on the object $\ovl M_U$, see Section \ref{section the geometric model -1 cy}. The observation below will be useful for describing the irreducible morphisms of $\rep(\Zc_m,\ka)$, see Proposition \ref{proposition irreducible morphism in rep}.

\begin{observation}\label{observation intervals}
Let $U = (u_1,u_2+2h\pi]\in\I$. From Figure \ref{figure ar quiver rep}, it is straightforward to see that the following statements hold.
\begin{enumerate}
\item $U_1\in\I$ if and only if $u_1^+\leq u_2+2h\pi\leq u_1+2\pi$, or equivalently $U\neq (u_1,u_1^++2\pi]$.
\item $U_2\in\I$ if and only if $U_2\neq\emp$, or equivalently $U\neq (u_1,u_1^+]$.
\item $\S^{-1}U\in\I$ if and only if $U\neq (u_1,u_1^++2\pi]$. The same holds for $\S U$.
\item $U^-\in\I$. Moreover, $U^- = (u_2^-,u_2]$ if and only if $U = (u_1,u_1^+]$, and $U^- = (u_2^-,u_2+2\pi]$ if and only if $U = (u_1,u_1^++2\pi]$.
\item $U'\in\I$. Moreover, $U' = (u_2^-,u_2]$ if and only if $U = (u_1,u_1^++2\pi]$, and $U' = (u_2^-,u_2+2\pi]$ if and only if $U = (u_1,u_1^++2\pi]$.
\end{enumerate}
\end{observation}

\section{Morphisms and factorisation properties}\label{section morphisms}

We describe the $\Hom$-spaces and the factorisation properties of the morphisms of $\rep(\Zc_m,\ka)$.

\subsection{Hom-hammocks}

We define the $\Hom$-hammocks and prove the following result.

\begin{proposition}\label{proposition hom spaces}
Let $U = (u_1,u_2+2h\pi], V = (v_1,v_2+2k\pi]\in\I$. Then 
\[
\Hom_{S^1}(\ovl{M}_U,\ovl{M}_V)\cong
\begin{cases}
\K^2 & \text{if $U = V = (u_1,u_1^++2\pi]$,}\\
\K & \text{if either $V\cap_L U\neq \emp$ or $(V-2\pi)\cap_L U\neq \emp$,}\\
0 & \text{otherwise.}
\end{cases}
\]
\end{proposition}

\begin{figure}[h]
\centering
\includegraphics[height = 5cm]{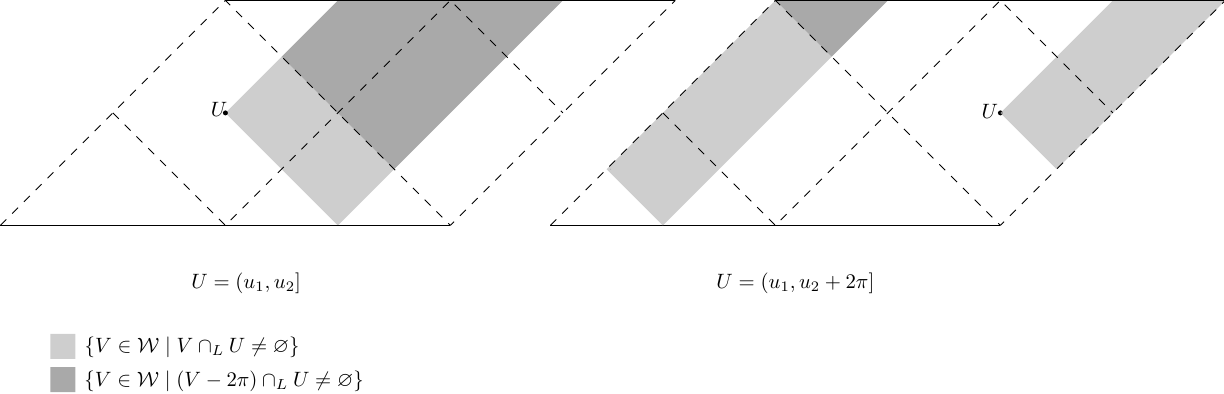}
\caption{For the interval $U\in\I$, the grey areas indicate the sets of intervals $V\in\I$ such that $V\cap_L U\neq\emp$ or $(V-2\pi)\cap_L U\neq\emp$.}
\label{figure hammocks rep 1}
\end{figure}

We need the following lemma, which relates the $\Hom$-spaces of $\rep(\Zc_m,\ka)$ to the intersections of intervals of $\I$. We refer to Figure \ref{figure hammocks rep 1} for an illustration in the coordinate system of $\rep(\Zc_m,\ka)$.

\begin{lemma}\label{lemma dimension hom spaces}
Let $U = (u_1,u_2+2h\pi], V = (v_1,v_2+2k\pi]\in\I$. The following statements hold.
\begin{enumerate}
\item If $(V+2n\pi)\cap_L U\neq \emp$ for some $n\in\Z$, then $n\in\{0,-1\}$. Therefore $\Hom_{S^1}(\ovl{M}_U,\ovl{M}_V)$ is at most two dimensional.	
\item We have that $V\cap_L U\neq \emp$ if and only if $v_1\leq u_1$ and $u_1^+\leq v_2+2k\pi\leq u_2+2h\pi$.	
\item We have that $(V-2\pi)\cap_L U\neq \emp$ if and only if $k = 1$ and $u_1^+\leq v_2\leq u_2+2h\pi$.
\item We have that $V\cap_L U\neq \emp$ and $(V-2\pi)\cap_L U\neq \emp$ if and only if $U = V = (u_1,u_1^++2\pi]$.
\end{enumerate}
\end{lemma}
\begin{proof}
We prove statement (1). If $n = 1$ then, since $0<u_1<2\pi$ and $2\pi<v_1+2\pi<4\pi$, we have that $(V+2\pi)\cap_L U = \emp$. Moreover, if $n\geq 2$ or $n\leq -2$, then $(V+2n\pi)\cap U = \emp$ and as a consequence $(V+2n\pi)\cap_L U = \emp$. Thus, $n\in\{-1,0\}$.

Statements (2) and (3) are straightforward, we prove (4). It is straightforward to check that if $U = V = (u_1,u_1^++2\pi]$ then $V\cap_L U\neq \emp$ and $(V-2\pi)\cap_L U\neq\emp$. If $V\cap_L U\neq \emp$ and $(V-2\pi)\cap_L U\neq\emp$, then $k = 1$, $u_1^+\leq v_2\leq u_2+2h\pi$, and $u_1^+\leq v_2+2\pi\leq u_2+2h\pi$. Thus, $h = 1$ and then $v_2\leq u_2$. Since $u_2+2h\pi = u_2+2\pi\leq u_1^++2\pi$, we have that $u_1^+\leq v_2\leq u_2\leq u_1^+$, i.e. $v_2 = u_1^+$. Note that $v_1\geq v_2^-+2k\pi-2\pi = u_1$ and, since $V\cap_L U\neq \emp$, also $v_1\leq u_1$. Therefore, $v_1 = u_1$. Moreover, we have that $u_2+2h\pi = u_2+2\pi\leq u_1^++2\pi$, i.e. $u_2\leq u_1^+$. Since $u_2+2\pi\leq u_1^++2\pi = v_2+2\pi\leq u_2+2\pi$, we obtain that $u_2 = u_1^+$. We conclude that $U = V = (u_1,u_1^++2\pi]$.
\end{proof}

Now we can prove Proposition \ref{proposition hom spaces}.

\begin{proof}[Proof of Proposition \ref{proposition hom spaces}]
By Theorem \ref{theorem indecomposable strings} we have that $\Hom_{S^1}(\ovl{M}_U,\ovl{M}_V)\cong \K^n$ with $n = |\{l\in\Z\mid (V+2l\pi)\cap_L U \neq \emp\}|$. Then the claim follows from Lemma \ref{lemma dimension hom spaces}.
\end{proof}

\begin{figure}[h]
\centering
\includegraphics[height = 5cm]{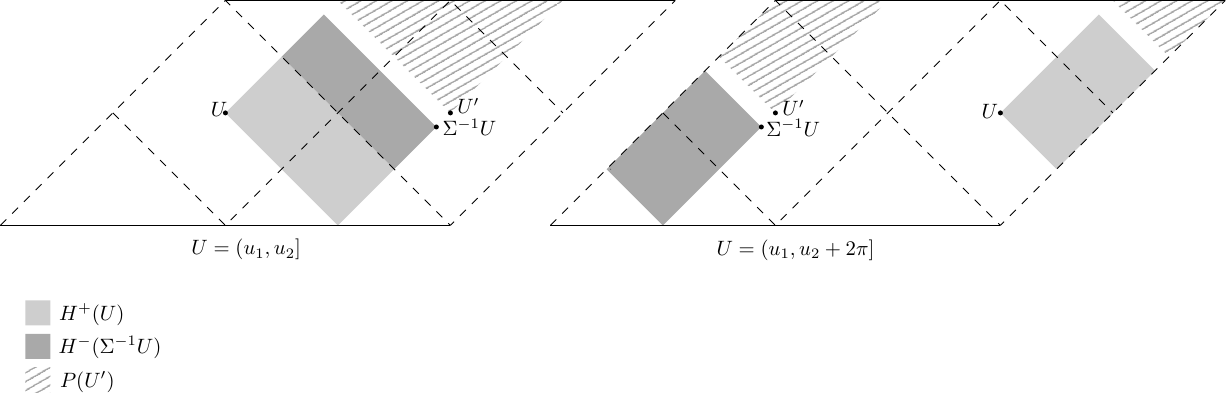}
\caption{The interval $U\in\I$ and its $\Hom$-hammocks, cf. Figure \ref{figure hammocks rep 1}.}
\label{figure hammocks rep 2}
\end{figure}

We define the $\Hom$-hammocks in $\rep(\Zc_m,\ka)$, Figure \ref{figure hammocks rep 2} provides an illustration.

\begin{definition}\label{definition hom hammocks}
Let $U = (u_1,u_2+2h\pi]\in\I$. We define the following sets.
\begin{align*}
H^+(U) & = 
\begin{cases}
\{(v_1,v_2]\in\I\mid v_1\leq u_1\text{ and }u_1^+\leq v_2\leq u_2\} & \text{if $h = 0$,}\\
\{(v_1,v_2+2\pi]\in\I\mid u_2\leq v_1\leq u_1 \text{ and } v_2\leq u_2\} & \text{if $h = 1$.}
\end{cases}\\[0.5cm]
H^-(U) & = 
\begin{cases}
\{(v_1,v_2]\in\I\mid u_1\leq v_1\leq u_2^-\text{ and } v_2\geq u_2\} & \text{if $h = 0$,}\\
\{(v_1,v_2+2\pi]\in\I\mid v_1\geq u_1 \text{ and } u_2\leq v_2\leq u_1\} & \text{if $h = 1$.}
\end{cases}\\[0.5cm]
P(U) & = 
\begin{cases}
\left\{(v_1,v_2+2k\pi]\in\I\middle| \parbox{6.2cm}{\centering $v_1-2k\pi\leq u_1$ and $v_2\geq u_2$, or $v_1\leq u_1$ and $v_2+2k\pi\geq u_2$}\right\} & \text{if $h = 0$,}\\[0.35cm]
\{(v_1,v_2+2\pi]\in\I\mid v_1\leq u_1 \text{ and } v_2\geq u_2\} & \text{if $h = 1$.}
\end{cases}
\end{align*}
We extend these definitions to $U = \emp$ by imposing $H^+(\emp) = H^-(\emp) = P(\emp) = \emp$.
\end{definition}

\begin{observation}\label{observation h+ h- p}
Let $U = (u_1,u_2+2h\pi]\in\I$ and $V = (v_1,v_2+2k\pi]\in\I$ be such that $V\in H^+(U)\sqcup H^-(\S^{-1}U)\sqcup P(U')$. The following statements hold, see Figure \ref{figure hammocks rep 2} for an illustration.
\begin{enumerate}
\item The sets $H^+(U)$, $H^-(\S^{-1}U)$, and $P(U')$ are pairwise disjoint.
\item $V\in H^+(U)$ if and only if $V\notin P(U')$ and $k = h$. Moreover, if $V\in H^+(U)$ then $V\cap_L U\neq\emp$.
\item $V\in H^-(\S^{-1}U)$ if and only if $V\notin P(U')$ and $k = 1-h$. Moreover, if $V\in H^-(\S^{-1}U)$ then $(V-2k\pi)\cap_L U\neq\emp$.
\item If $U = (u_1,u_1^+]$, then $P(U') = \{(u_1,u_1^++2\pi]\}$.
\item If $U = (u_1,u_1^++2\pi]$, then $H^+(U) = H^-(\S^{-1}U) = \emp$.
\end{enumerate}
\end{observation}

Now we can prove the following result.

\begin{proposition}\label{proposition hom hammocks}
Let $U, V\in\I$, then $\Hom_{S^1}(\ovl{M}_U,\ovl{M}_V)\neq 0$ if and only if $V\in H^+(U)\sqcup H^-(\S^{-1}U)\sqcup P(U')$.
\end{proposition}
\begin{proof}
We denote $U = (u_1,u_2+2h\pi]$ and we assume that $h = 0$, if $h = 1$ the proof is analogous. It is straightforward to check that the following equality holds.
\begin{align*}
H^+(U)\sqcup H^-(\S^{-1}U) \sqcup P(U') = & \{(v_1,v_2]\in\I\mid u_1^+\leq v_2\leq u_2 \text{ and } v_1\leq u_1\}\sqcup\\
& \{(v_1,v_2+2\pi]\in\I\mid u_1^+\leq v_2\leq u_2\}
\end{align*}
Let $V\in\I$. By Lemma \ref{lemma dimension hom spaces} and Proposition \ref{proposition hom spaces}, we know that $\Hom_{S^1}(\ovl{M}_U,\ovl{M}_V)\neq 0$ if and only if $V\cap_L U\neq \emp$ or $(V-2\pi)\cap_L U\neq \emp$. Thus, it is straightforward to check that
\begin{align*}
\{V\in\I\mid \Hom_{S^1}(\ovl{M}_U,\ovl{M}_V)\neq \emp\} = & \{(v_1,v_2]\in\I\mid v_1\leq u_1\text{ and }u_1^+\leq v_2\leq u_2\}\sqcup\\
& \{(v_1,v_2+2\pi]\in\I\mid u_1^+\leq v_2\leq u_2\}.
\end{align*}
We conclude that $V\in\I$ is such that $\Hom_{S^1}(\ovl{M}_U,\ovl{M}_V)\neq \emp$ if and only if $V\in H^+(U)\sqcup H^-(\S^{-1}U)\sqcup P(U')$.
\end{proof}

With Proposition \ref{proposition factorisation through projectives} we will prove that $P(U')$ determines the set of intervals $V\in\I$ such that $\Proj(\ovl M_U,\ovl M_V) \neq 0$.

\subsection{Factorisation properties}\label{section factorisation properties rep}

We prove the factorisation properties of the morphisms of $\rep(\Zc_m,\ka)$. After stabilising, the proposition below can be reformulated as in Proposition \ref{proposition factorisation properties}.

\begin{proposition}\label{proposition factorisation properties in rep}
Let $U,V,W\in\I$ be such that $\ovl{M}_U$, $\ovl{M}_V$, and $\ovl{M}_W$ are non-isomorphic. Assume that there exist non-zero morphisms $f\colon \ovl{M}_U\to \ovl{M}_V$ and $g\colon \ovl{M}_V\to \ovl{M}_W$. Assume that one of the following conditions holds.
\begin{enumerate}
\item $V\cap_L U\neq \emp$, $W\cap_L V\neq\emp$, and $W\cap_L U\neq\emp$.
\item $V\cap_L U\neq \emp$, $(W-2\pi)\cap_L V\neq\emp$, and $(W-2\pi)\cap_L U\neq\emp$.
\item $(V-2\pi)\cap_L U\neq \emp$, $W\cap_L V\neq\emp$, and $(W-2\pi)\cap_L U\neq\emp$.
\end{enumerate}
Then $gf\neq 0$.
\end{proposition}
\begin{proof}
By assumption, we have that $(V+2n\pi)\cap_L U\neq\emp$, $(W+2l\pi)\cap_L (V+2n\pi)\neq\emp$, and $(W+2l\pi)\cap_L U\neq\emp$ for some $l,n\in\Z$. Since $U$, $V$, and $W$ are pairwise non-isomorphic, by Proposition \ref{proposition hom spaces} we have that such $l$ and $n$ are unique. Let $\phi\colon M_U\to M_{V+2n\pi}$ and $\psi\colon M_{V+2n\pi}\to M_{W+2l\pi}$ be standard morphisms, see Definition \ref{definition standard morphisms}. Since $(W-2l\pi)\cap_L U\neq \emp$, $\psi\phi\neq 0$, and then $\ovl\psi\hspace{0.05 cm}\ovl\phi = \ovl{\psi\phi}\neq 0$, see Lemma \ref{lemma morphisms in rep}. Moreover, $f = \l\ovl\phi$ and $g = \mu\ovl\psi$ for some $\l,\mu\in\K^*$, and as a consequence $gf = \l\mu\ovl{\psi}\hspace{0.05 cm}\ovl{\phi}\neq 0$. We conclude that $gf\neq 0$.
\end{proof}

\begin{proposition}\label{proposition zero compositions in rep}
Let $U$, $V$, $W$, $f$, $g$ be as in Proposition \ref{proposition factorisation properties in rep}. Assume that one of the following conditions holds.
\begin{enumerate}
\item $(V-2\pi)\cap_L U\neq \emp$ and $(W-2\pi)\cap_L V\neq\emp$.
\item $V\cap_L U\neq \emp$, $(W-2\pi)\cap_L V\neq\emp$, and $W\cap_L U\neq\emp$.
\item $(V-2\pi)\cap_L U\neq \emp$, $W\cap_L V\neq \emp$, and $W\cap_L U\neq \emp$.
\end{enumerate}
Then $gf = 0$.
\end{proposition}
\begin{proof}
We proceeding similarly as in the argument of Proposition \ref{proposition factorisation properties in rep}. By assumption we have that $(V+2n\pi)\cap_L U\neq \emp$ and $(W+2l\pi)\cap_L (V+2n\pi)\neq \emp$ for some $l,n\in\Z$. Let $\phi\colon M_U\to M_{V+2n\pi}$ and $\psi\colon M_{V+2n\pi}\to M_{W+2l\pi}$ be standard morphisms, we obtain that $gf = \l\ovl{\psi\phi}$ for some $\l\in\K^*$. Note that $(W-2l\pi)\cap_L U = \emp$. Indeed, if condition (1) holds, we have that $l = -2$ and then $(W+2l\pi)\cap_L U = (W-4\pi)\cap_L U = \emp$. If condition (2) or (3) holds, then $l = -1$ and, since $W\cap_L U\neq \emp$, we have that $(W-2l\pi)\cap_L U = (W-2\pi)\cap_L U = \emp$. Indeed, if both $W\cap_L U\neq\emp$ and $(W-2\pi)\cap_L U\neq\emp$, then, by Lemma \ref{lemma dimension hom spaces}, $U = W = (u_1,u_1^++2\pi]$ and then $\ovl M_U\cong\ovl M_V$, giving a contradiction. As a consequence, $\psi\phi = 0$ and we conclude that $gf = \l\ovl{\psi\phi} = 0$.
\end{proof}

The following lemma describes when there is a monomorphism or an epimorphism between two indecomposable objects of $\rep(\Zc_m,\ka)$. It will be useful for proving Propositions \ref{proposition enough injective objects} and \ref{proposition simple objects}.

\begin{lemma}\label{lemma mono and epi}
Let $U = (u_1,u_2+2h\pi],V = (v_1,v_2+2k\pi]\in\I$. The following statements hold.
\begin{enumerate}
\item Assume that $U\neq (u_1,u_1^++2\pi]$. Then there exists a monomorphism $\ovl M_U\to \ovl M_V$ if and only if there exists $l\in\{0,-1\}$ such that $(V+2l\pi)\cap_L U\neq\emp$ and $u_2+2h\pi = v_2+2(k+l)\pi$.
\item Assume that $V\neq (v_1,v_1^++2\pi]$. Then there exists an epimorphism $\ovl M_U\to \ovl M_V$ if and only if $V\cap_L U\neq\emp$ and $u_1 = v_1$.
\end{enumerate} 
\end{lemma}
\begin{proof}
We prove statement (1). Assume that $(V+2l\pi)\cap_L U\neq\emp$ and $u_2+2h\pi = v_2+2(k+l)\pi$ for a unique $l\in\{0,-1\}$, then $\Hom_{S^1}(\ovl M_U,\ovl M_V)\cong\K$. Let $f\colon \ovl M_U\to\ovl M_V$ be non-zero, we prove that $f$ is a monomorphism. Consider a non-zero morphism $g\colon \ovl M_W\to \ovl M_U$ with $W\in\I$, we check that $fg \neq 0$. Since $\Hom (\ovl M_W,\ovl M_U)\cong\K$, we have that either $U\cap_L W\neq\emp$ or $(U-2\pi)\cap_L W\neq\emp$. If $(U-2\pi)\cap_L W\neq\emp$ and $l = -1$, then $h = 1$ and we obtain that $u_2+2\pi = u_2+2h\pi = v_2+2(k+l)\pi = v_2+2k\pi-2\pi$, i.e. $u_2+4\pi = v_2+2k\pi$, which is impossible because $k\in\{0,1\}$. 
	
Thus, we have the following possibilities: $U\cap_L W\neq\emp$ and $l = 0$, $U\cap_L W\neq\emp$ and $l = -1$, or $(U-2\pi)\cap_L W\neq\emp$ and $l = 0$. Since $u_2+2h\pi = v_2+2(k+l)\pi$, in the first case we obtain that $V\cap_L W\neq\emp$, and for the remaining cases $(V-2\pi)\cap_L W\neq\emp$, see Figure \ref{figure monomorphism}. By Proposition \ref{proposition factorisation properties in rep} we obtain that $fg\neq 0$. This proves that $f$ is a monomorphism.
	
\begin{figure}[h]
\centering
\includegraphics[height = 1.2cm]{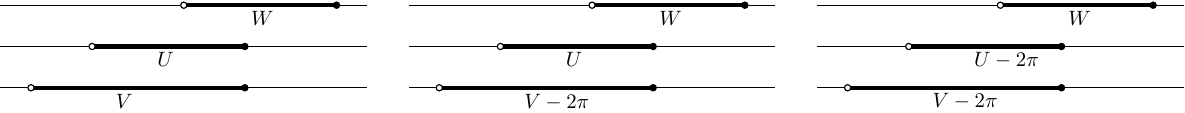}
\caption{The intersections of the argument of Lemma \ref{lemma mono and epi}.}
\label{figure monomorphism}
\end{figure}
	
Now assume that there exists a monomorphism $f\colon\ovl M_U\to \ovl M_V$. Since $f\neq 0$, there exists a unique $l\in\{0,-1\}$ such that $(V+2l\pi)\cap_L U\neq\emp$, we prove that $u_2+2h\pi = v_2+2(k+l)\pi$. There exists a monomorphism $g\colon \ovl M_{(u_2^-,u_2]}\to \ovl M_U$. Indeed, it is straightforward to check that $(U-2h\pi)\cap_L (u_2^-,u_2]\neq\emp$ and $u_2 = u_2+2(h-h)\pi$. Thus, $fg\colon \ovl M_{(u_2^-,u_2]}\to \ovl M_V$ is a monomorphism, and as a consequence $\Hom_{S^1}(\ovl M_{(u_2^-,u_2]},\ovl M_V)\cong\K$. We have that either $V\cap_L (u_2^-,u_2]\neq\emp$ or $(V-2\pi)\cap_L (u_2^-,u_2]\neq\emp$. For both cases, by Lemma \ref{lemma dimension hom spaces}, we have that $v_2 = u_2$. Moreover, since $(V+2l\pi)\cap_L U\neq\emp$, we have that $u_1^+\leq v_2+2(k+l)\pi\leq u_2+2h\pi$. Thus, $v_2+2(k+l)\pi = u_2+2(k+l)\pi\in U$. Since $U \neq (u_1,u_1^++2\pi]$, $U$ contains only one element, namely $u_2+2h\pi$, which is equal to $u_2+2h\pi$ up to a multiple of $2\pi$. Therefore, $v_2+2(k+l)\pi = u_2+2h\pi$. This concludes the argument of (1). 
	
Now we prove statement (2). Similarly as above, we can prove that if $V\cap_L U\neq\emp$ and $v_1 = u_1$, then $\Hom_{S^1}(\ovl M_U,\ovl M_V)\cong\K$ and any non-zero morphism $\ovl M_U\to\ovl M_V$ is an epimorphism. Now assume that there exists an epimorphism $f\colon \ovl M_U\to\ovl M_V$, we prove that $V\cap_L U\neq\emp$ and $v_1 = u_1$. Since $(v_1,v_1^+]\cap_L V\neq\emp$, there is an epimorphism $g\colon \ovl M_V\to\ovl M_{(v_1,v_1^+]}$. Thus, $gf\colon \ovl M_U\to \ovl M_{(v_1,v_1^+]}$ is an epimorphism, and in particular $gf\neq 0$. As a consequence, either $(v_1,v_1^+]\cap_L U\neq\emp$ or $((v_1,v_1^+]-2\pi)\cap_L U\neq\emp$. Note that the second case is impossible because $(v_1,v_1^+]\subseteq(0,2\pi)$. We obtain that $V\cap_L U\neq\emp$, indeed, by Proposition \ref{proposition zero compositions in rep}, if $(V-2\pi)\cap_L U\neq\emp$ then $gf = 0$, giving a contradiction. Moreover, since $(v_1,v_1^+]\cap_L U\neq\emp$, we have that $v_1\leq u_1$ and $v_1^+\geq u_1^+$, i.e. $v_1 = u_1$. This concludes the argument of (2). 
\end{proof}

\section{Projective-injective objects}\label{section projective-injective objects rep}

By Theorems \ref{theorem rep is wide} and \ref{theorem decomposition theorem for rep}, $\rep(\Zc_m,\ka)$ is a Krull--Schmidt wide subcategory of $\Rep(S^1,\ka)$. We now prove that $\rep(\Zc_m,\ka)$ is a Frobenius subcategory of $\Rep^{\pwf}(S^1,\ka)$ and that for each $z\in S^1$ we have that $P_z = I_z$, where $P_z$ and $I_z$ are respectively the indecomposable projective and indecomposable injective objects at $z$. We refer to Section \ref{section frobenius categories} for some background on Frobenius categories. 

Given $z\in\Zc_m$, we denote the object $\ovl M_{(z,z^++2\pi]}$ by $P_z$ or $I_z$. The following theorem is analogous to \cite[Remark 6.5(2)]{RZ} and follows directly from Propositions \ref{proposition projective injective objects} and \ref{proposition enough injective objects}.

\begin{theorem}\label{theorem frobenius}
The category $\rep(\Zc_m,\ka)$ is Frobenius and its indecomposable projective-injective objects are exactly those isomorphic to $P_z = I_z$ for some $z\in\Zc_m$.
\end{theorem}

\begin{proposition}\label{proposition projective injective objects}
For each $z\in\Zc_m$ the object $P_z = I_z$ is projective and injective in $\rep(\Zc_m,\ka)$.
\end{proposition}
\begin{proof}
Let $z\in\Zc_m$, we prove that $P_z$ is projective. The proof that $I_z$ is injective is dual. Consider a short exact sequence $0\lora L\ovs{f}{\lora} M\ovs{g}{\lora} P_z\lora 0$ of $\rep(\Zc_m,\ka)$, we show that it splits, i.e. $g$ is a split epimorphism. We divide the proof into claims.
	
\emph{Claim 1.} There exists a direct summand of $M$ isomorphic to $P_z$.
	
Consider the following commutative diagram of vector spaces with split exact rows.
\[
\begin{tikzcd}
0 \ar[r]& L(z^+)\ar[r,"f(z^+)"]\ar[d,"L(\om_{z^+})"] & M(z^+)\ar[r,"g(z^+)"]\ar[d,"M(\om_{z^+})"] & P_z(z^+)\ar[r]\ar[d,"P_z(\om_{z^+})"] & 0	\\
0 \ar[r]& L(z^+)\ar[r,"f(z^+)"] & M(z^+)\ar[r,"g(z^+)"] & P_z(z^+)\ar[r] & 0
\end{tikzcd}
\]
We recall that, by Definition \ref{definition string}, $P_z(z^+) = \K b_{1,z^+}\oplus \K b_{2,z^+}$,  $P_z(\om_{z^+})(b_{1,z^+}) = b_{2,z^+}$, and $P_z(\om_{z^+})(b_{2,z^+}) = 0$. Assume that $M(\om_{z^+}) = 0$. Then $P_z(\om_{z^+})g(z^+) = 0$, and this contradicts the fact that $g(z^+)$ is a split epimorphism and $P_z(\om_{z^+})\neq 0$. Therefore $M(\om_{z^+})\neq 0$. Now consider the decomposition of $M$ into indecomposable direct summands $M = \bigoplus_{i = 1}^n M_i$, we can write $M(\om_{z^+}) = \bigoplus_{i = 1}^n M_i(\om_{z^+})$. Since $M(\om_{z^+}) \neq 0$, there exists $i\in\{1,\dots,n\}$ such that $M_i(\om_{z^+})\neq 0$, i.e. $M_i\cong P_z$. This concludes the argument of Claim 1.
	
Now we introduce some notation. After reordering the summands of $M$, we can write $M = \bigoplus_{i = 1}^k M_i \oplus \bigoplus_{j = k+1}^n M_j$ where $M_i = P_{z}$ for each $i\in\{1,\dots,k\}$ and $M_j\not\cong P_{z}$ for each $j\in\{k+1,\dots,n\}$. We also write $g\colon M\to P_{z}$ as $g = \left(\begin{smallmatrix}g_1 & \dots & g_k & g_{k+1} & \dots & g_n\end{smallmatrix}\right)$. We show that there exists $i\in\{1,\dots,k\}$ such that $g_i\colon M_i\to P_{z}$ is an isomorphism. First, we need to check the following fact. 
	
\emph{Claim 2.} For each $i\in\{1,\dots,k\}$ and $j\in\{k+1,\dots,n\}$, we have that $g_i(z^+) = \left(\begin{smallmatrix}\a_i & 0\\ \b_i & \a_i\end{smallmatrix}\right)$ and  $g_j(z^+) = \left(\begin{smallmatrix} 0 & \g_j\end{smallmatrix}\right)^T$ for some $\a_i,\b_i,\g_j\in \K$.

We recall that for each $i\in\{1,\dots,k\}$ and $j\in\{k+1,\dots,n\}$ the following diagrams commute.
\[
\begin{matrix}
\begin{tikzcd}
M_i(z^+) \ar[r,"g_i(z^+)"]\ar[d,"M_i(\om_{z^+})"]& P_{z}(z^+)\ar[d,"P_z(\om_{z^+})"]\\	
M_i(z^+) \ar[r,"g_i(z^+)"]& P_{z}(z^+)
\end{tikzcd}
&
\begin{tikzcd}
M_j(z^+) \ar[r,"g_j(z^+)"]\ar[d,"M_j(\om_{z^+})"]& P_{z}(z^+)\ar[d,"P_z(\om_{z^+})"]\\
M_j(z^+) \ar[r,"g_j(z^+)"]& P_{z}(z^+)
\end{tikzcd}
\end{matrix}
\]
Then the linear map $g_i(z^+)\colon M_i(z^+) = \K b_{1,z^+}\oplus \K b_{2,z^+}\to P_z(z^+) = \K b_{1,z^+}\oplus \K b_{2,z^+}$ is of the form
$g_i(z^+) = \left(\begin{smallmatrix}\a_i & 0\\ \b_i & \a_i\end{smallmatrix}\right)$ for some $\a_i,\b_i\in\K$. Moreover, since $M_j\not\cong P_z$, we have that $M_j(z^+) = 0$ or $M_j(z^+)\cong \K$ and in both cases $M_j(\om_{z^+}) = 0$. Therefore, $g_j(z^+) = \left(\begin{smallmatrix} 0 & \g\end{smallmatrix}\right)^T$ for some $\g\in \K$. 

\emph{Claim 3.} There exists $i\in\{1,\dots,k\}$ such that $g_i(z^+)\colon M_i(z^+)\to P_z(z^+)$ is an isomorphism.

We show that $\a_i\neq 0$ for some $i\in\{1,\dots,k\}$. Indeed, if $\a_i = 0$ for all $i\in\{1,\dots,k\}$, then 
\[
g(z^+) = 
\left(
\begin{matrix} 
0 & 0 & \cdots & 0 	 & 0 & 0 & \cdots & 0\\
\b_1 & 0 & \cdots & \b_k  & 0 & \g_{k+1} & \cdots & \g_n
\end{matrix}
\right)
\]
and this contradicts the fact that $g(z^+)$ is a split epimorphism. As a consequence, there exists $i\in\{1,\dots,k\}$ such that $\a_i\neq 0$, and therefore $g_i(z^+)$ is an isomorphism.

\emph{Claim 4.} Let $i\in\{1,\dots,k\}$ be such that $g_i(z^+)$ is an isomorphism. We have that $g_i\colon M_i\to P_z$ is an isomorphism.

We prove that $g_i(x)\colon M_i(x)\to P_z(x)$ is an isomorphism for each $x\in S^1$. Let $x\in S^1$ and consider the following commutative diagram.
\[
\begin{tikzcd}
M_i(x)\ar[r,"g_i(x)"]\ar[d,"M_i(g_{xz^+})"] & P_z(x)\ar[d,"P_z(g_{xz^+})"]\\
M_i(z^+)\ar[r,"g_i(z^+)"] & P_z(z^+)
\end{tikzcd}
\]
If $z<x\leq z^+$ are in cyclic order then, since $M_i(g_{xz^+})$ and $P_z(g_{xz^+})$ are isomorphisms by Definition \ref{definition string}, and since $g_i(z^+)$ is an isomorphism by Claim 2, we obtain that $g_i(x)$ is an isomorphism. Now assume that $z^+<x\leq z$ are in cyclic order. We recall that $M_i(x) = \K b_{1,x} = P_z(x)$ and $M_i(g_{xz^+})(b_{1,x}) = b_{1,z^+} = P_z(g_{xz^+})(b_{1,x})$. Since the diagram above commutes, we obtain that $g_i(x)\neq 0$, and therefore $g_i(x)$ is an isomorphism. Thus, $g_i(x)$ is an isomorphism for each $x\in S^1$, i.e. $g_i$ is an isomorphism.

\emph{Claim 5.} The morphism $g\colon M\to P_z$ is a split epimorphism.

From Claim 4 we know that $g_i\colon M_i\to P_z$ is an isomorphism in $\rep(\Zc_m,\ka)$. We conclude that $g\colon M\to P_z$ is a split epimorphism, i.e. the short exact sequence $0\lora L\ovs{f}{\lora} M\ovs{g}{\lora} P_z\lora 0$ splits. 
\end{proof} 

The following are important properties for the projective and injective objects of $\rep(\Zc_m,\ka)$.

\begin{proposition}\label{proposition enough injective objects}
The category $\rep(\Zc_m,\ka)$ has enough projectives and enough injectives. Moreover, each indecomposable projective or injective object of $\rep(\Zc_m,\ka)$ is isomorphic to $P_z = I_z$ for some $z\in\Zc_m$.
\end{proposition}
\begin{proof}
By Theorem \ref{theorem decomposition theorem for rep}, each object of $\rep(\Zc_m,\ka)$ decomposes as a finite direct sum of indecomposable objects, thus it is enough to prove that each indecomposable object of $\rep(\Zc_m,\ka)$ has a projective cover and an injective envelope. Let $U = (u_1,u_2+2h\pi]\in\I$, by Lemma \ref{lemma mono and epi} there exist an epimorphism $P_{u_1}\to \ovl{M}_U\to 0$ and a monomorphism $0\to I_{u_2^-}\to \ovl{M}_U$. 

Now we prove that the indecomposable projective or injective objects are all of the form $P_z = I_z$ for some $z\in\Zc_m$. Let $U = (u_1,u_2+2h\pi]\in\I$ and $\ovl{M}_{U}$ be projective in $\rep(\Zc_m,\ka)$, then there exists an epimorphism $P_{u_1}\to \ovl{M}_{U}\to 0$. Since $\ovl{M}_{U}$ is projective, this epimorphism splits and then $\ovl{M}_{U}$ is a direct summand of $P_{u_1}$. Therefore, $\ovl{M}_{U}\cong P_{u_1}$. If $\ovl M_U$ is injective, we can proceed dually.
\end{proof}

With Proposition \ref{proposition hom hammocks} we proved that, for any pair of intervals $U,V\in\I$, $\Hom_{S^1}(\ovl M_U,\ovl M_V)\neq 0$ if and only if $V\in H^+(U)\sqcup H^-(\S^{-1}U)\sqcup P(U')$. Now we characterise the intervals $V\in\I$ for which $\Proj(\ovl M_U,\ovl M_V)\neq 0$.

\begin{proposition}\label{proposition factorisation through projectives}
Let $U,V\in\I$. We have that $\Proj(\ovl{M}_U,\ovl{M}_V)\neq 0$ if and only if $V\in P(U')$.
\end{proposition}
\begin{proof}
We denote $U = (u_1,u_2+2h\pi]$ and $V = (v_1,v_2+2k\pi]$. If $\ovl M_U$ is projective, then $\Proj(\ovl M_U,\ovl M_V) = \Hom_{S^1}(\ovl M_U,\ovl M_V)\neq 0$ if and only if $V\in P(U')$, see Observation \ref{observation h+ h- p}. Now assume that $\ovl M_U$ is not projective. Let $f\colon \ovl M_U\to \ovl M_V$ be a non-zero morphism such that $f = hg$ for some $g\colon \ovl M_U\to P$ and $h\colon P\to \ovl M_V$, where $P$ is a projective object of $\rep(\Zc_m,\ka)$. Since $\ovl M_U$ is not projective, by Proposition \ref{proposition hom spaces}, $\Hom_{S^1}(\ovl M_U,\ovl M_V)\cong\K$. Therefore, there exists an indecomposable direct summand $P_z$ of $P$ such that $f = \b\a$ for some non-zero morphisms $\a\colon \ovl M_U\to P_z$ and $\b\colon P_z\to \ovl M_U$.  

If $h = 0$, since $(z,z^++2\pi]\in P(U')$, we have that $u_1^+\leq z^+\leq u_2$ and $\left((z,z^++2\pi]-2\pi\right)\cap_L U\neq \emp$. We have that $V\cap_L (z,z^++2\pi]\neq\emp$ and $(V-2\pi)\cap_L U\neq\emp$, otherwise, by Proposition \ref{proposition zero compositions in rep}, $\b\a = 0$. Thus, $v_1\leq z\leq  u_2^-$ and $v_2\geq u_1^+$, i.e. $V\in P(U')$.  

If $h = 1$, we have the following possibilites: either $(V-2\pi)\cap_L U\neq\emp$ or $V\cap_L U\neq\emp$. In the first case, $k = 1$ and $v_2\geq u_1^+$, and as a consequence $V\in P(U')$. In the second case, since $V\cap_L U\neq\emp$, $(z,z^++2\pi]\cap_L U\neq\emp$ and $V\cap (z,z^++2\pi]\neq\emp$, otherwise $\b\a = 0$ by Proposition \ref{proposition zero compositions in rep}. Therefore, $v_2+2k\pi\geq u_1^+$ and, since $z^++2\pi\leq u_2+2\pi$, $v_1\leq z\leq u_2^-$, and then $V\in P(U')$. This proves that if $\Proj(\ovl M_U,\ovl M_U)\neq 0$ then $V\in P(U')$.

Now we prove that if $V\in P(U')$ then $\Proj(\ovl M_U,\ovl M_V)\neq 0$. Let $f\colon \ovl M_U\to \ovl M_V$ be a non-zero morphism, we show that $f$ factors through a projective object. Assume that $h = 0$. Note that $\left((u_2^-,u_2+2\pi]-2\pi\right)\cap_L U \neq \emp$. Since $V\in P(U')$, we have that $v_1\leq u_2^-$ and $v_2\geq u_1^+$, and then $V\cap (u_2^-,u_2+2\pi]$ and $(V-2\pi)\cap_L U\neq \emp$. Thus, by Proposition \ref{proposition factorisation properties in rep}, $f$ factors through $P_{u_2^-}$. 

Now assume that $h = 1$. We have that $(u_2^-,u_2+2\pi]\cap U\neq \emp$. Since $V\in P(U')$, we have the following possibilities: $v_1\leq u_2^-$ and $v_2+2k\pi\geq u_1^+$, or $v_1-2k\pi\leq u_2^-$ and $v_2\geq u_1^+$. Thus, we have that $V\cap_L (u_2^-,u_2+2\pi]\neq\emp$ and $V\cap_L U\neq\emp$, or $(V-2\pi)\cap_L (u_2^-,u_2+2\pi]\neq\emp$ and $(V-2\pi)\cap_L U\neq\emp$. In both cases $f$ factors through $P_{u_2^-}$. We conclude that any non-zero morphism of $\Hom_{S^1}(\ovl M_U,\ovl M_V)$ factors through a projective object, i.e. $\Proj(\ovl M_U,\ovl M_V)\cong\K$.
\end{proof}

\section{Exact sequences}

We compute the middle terms of certain exact sequences of $\rep(\Zc_m,\ka)$ having indecomposable outer terms. This computation will be useful in Section \ref{section simple objects rep} for proving that $\rep(\Zc_m,\ka)$ is uniserial, and in Section \ref{section irreducible morphisms and almost split sequences in rep} for describing the almost split sequences. Given $U\in\I$, we refer to Definition \ref{definition intervals} for the notation $U^-$.

\begin{setup}\label{setup short exact sequences}
Let $U = (u_1,u_2+2h\pi], V = (v_1,v_2+2k\pi]\in\I$ be such that $\ovl M_U$ is not projective, $\Hom_{S^1}(\ovl M_{U^-},\ovl M_V)\cong\K$, and $\Proj(\ovl M_{U^-},\ovl M_V) = 0$. Since $\Hom_{S^1}(\ovl M_{U^-},\ovl M_V)\cong\K$, there exists a unique $l\in\{0,-1\}$ such that $(V+2l\pi)\cap_L U^-\neq\emp$. We denote $I = (v_1,u_2+2(h-l)\pi]$ and $J = (u_1,v_2+2(k+l)\pi]$.
\end{setup}

\begin{figure}[h]
\centering
\includegraphics[height = 2.5cm]{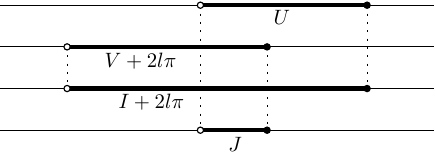}
\caption{The intervals $I$ and $J$ of Setup \ref{setup short exact sequences}.}
\label{figure intervals I J}
\end{figure}

The following lemmas will be useful for proving Proposition \ref{proposition short exact sequence}.

\begin{lemma}\label{lemma I J}
Keeping Setup \ref{setup short exact sequences}, the following statements hold.
\begin{enumerate}
\item We have that $I = (U-2l\pi)\cup V$ and $J = U\cap (V+2l\pi)$.
\item We have that $I\in\I$. Moreover, $\ovl M_I$ is projective if and only if $u_2+2(h-l)\pi = v_1^++2\pi$.
\item We have that $J\in\I$ if and only if $J\neq \emp$, or equivalently $v_2+2(k+l)\pi\neq u_1$. 
\item We have that $V = \S^{-1}(U^-)$ if and only if $J = \emp$ and $\ovl M_I$ is projective.
\item There exist a monomorphism $\ovl M_U\to \ovl M_I$ and an epimorphism $\ovl M_I\to \ovl M_V$. Moreover, if $J\neq\emp$, there exist an epimorphism $\ovl M_U\to \ovl M_J$ and a monomorphism $\ovl M_J\to \ovl M_V$.
\end{enumerate}
\end{lemma}
\begin{proof}
Statement (1) is straightforward and is illustrated by Figure \ref{figure intervals I J}. We prove statement (2). Since $V\subseteq I$, we have that $I\neq\emp$. To prove that $I\in\I$, we check that $u_2+2(h-l)\pi\leq v_1^++2\pi$. First we show that $h-l\in\{0,1\}$. Since $h,-l\geq 0$, we have that $h-l\geq 0$. Moreover, if $l = 0$ then $h-l = h\leq 1$, and if $l = -1$ then $h = 0$, otherwise $\Proj(\ovl M_{U^-},\ovl M_V)\neq 0$, see Figures \ref{figure hammocks rep 1} and \ref{figure hammocks rep 2}, giving a contradiction. Thus, $h-l\in\{0,1\}$.

Now, if $h-l = 0$, then $u_2+2(h-l)\pi = u_2\leq v_1^++2\pi$. If $h-l = 1$, we have the following possibilities: $h = 0$ and $l = -1$, or $h = 1$ and $l = 0$. We prove that in both cases $u_2\leq v_1^+$. If $h = 0$ and $l = -1$, we also have $k = 1$, and, since $(V+2l\pi)\cap_L U^-\neq\emp$, this implies that $v_2 = v_2+2(k+l)\pi\geq u_1$. If $h = 1$ and $l = 0$, we have that $v_2+2k\pi\geq u_1$. Since $\Proj(\ovl M_{U^-},\ovl M_V) = 0$, i.e. $V\notin P((U^-)')$ where $(U^-)' = (u_2^{--}, u_1+2(1-h)\pi]$, we have that $v_1\not\leq u_2^{--}$, see Definition \ref{definition hom hammocks}. Thus, $v_1\geq u_2^-$, i.e. $u_2\leq v_1^+$ and then $u_2+2(h-l)\pi = u_2+2\pi\leq v_1^++2\pi$. We obtain that $I\in\I$. Moreover, $\ovl M_I$ is projective if and only if $u_2+2(h-l)\pi = v_1^++2\pi$.

Now we prove statement (3). Since $(V+2l\pi)\cap_L U^-\neq\emp$, we have that $v_2+2(k+l)\pi\leq u_2^-+2h\pi\leq u_2+2h\pi\leq u_1^++2\pi$, and then $v_2+2(k+l)\pi\leq u_1^++2\pi$. Thus, either $J = \emp$, i.e. $v_2+2(k+l)\pi = u_1$, or $J\neq\emp$ and then $J\in\I$.

We prove statement (4). Assume that $V = (v_1,v_2+2k\pi] = \S^{-1}(U^-) = (u_2^-,u_1+2(1-h)\pi]$. Then $k = 1-h = -l$, $u_2^- = v_1$, i.e. $u_2 = v_1^+$, and $v_2 = u_1$. Therefore, $u_2+2(h-l)\pi = v_1^++2\pi$ and $v_2+2(k+l)\pi = u_1$. By statements (2) and (3) we have that $\ovl M_I$ is projective and $J = \emp$. Now assume that the converse holds, i.e. $u_1 = v_2+2(k+l)\pi$ and $u_2+2(h-l)\pi = v_1^++2\pi$. Then $v_1 = u_2^-$, $v_2 = u_1$, $k+l = 0$, and $h-l = 1$. As a consequence, $v_2+2k\pi = v_2+2(1-h)\pi = u_1+2(1-h)\pi$. We conclude that $V = \S^{-1}(U^-)$.

Finally, we prove statement (5). We check that $U$, $V$, $I$ and $J$ satisfy the conditions of Lemma \ref{lemma mono and epi}. By Figure \ref{figure intervals I J} it is straightforward to see that $(I+2l\pi)\cap_L U\neq\emp$ and $(V+2l\pi)\cap_L(I+2l\pi)\neq\emp$, i.e. $V\cap_L I\neq\emp$, and, if $J\neq\emp$, $J\cap_L U\neq\emp$ and $(V+2l\pi)\cap_L J\neq\emp$. Moreover, it is straightforward to check that the endpoints of $U$, $V$, $I$ and $J$ satisfy the remaining conditions of Lemma \ref{lemma mono and epi}. Thus, there exist monomorphisms $\ovl M_U\to \ovl M_I$ and $\ovl M_J\to \ovl M_V$, and epimorphisms $\ovl M_U\to \ovl M_J$ and $\ovl M_I\to \ovl M_V$. 
\end{proof}

\begin{lemma}\label{lemma short exact sequence}
Keeping Setup \ref{setup short exact sequences}, the following statements hold. 
\begin{enumerate}
\item If $J\neq\emp$, then $\dim\ovl M_U(x)+\dim\ovl M_V(x) =\dim\ovl M_I(x)+\dim\ovl M_J(x)$ for each $x\in S^1$.
\item If $J = \emp$, then $\dim\ovl M_U(x)+\dim\ovl M_V(x) =\dim\ovl M_I(x)$ for each $x\in S^1$.
\end{enumerate}
\end{lemma}
\begin{proof}
We prove statement (1). Let $x\in S^1$. Since $\ovl M_U$ and $\ovl M_V$ are not projective, we have that $\dim\ovl M_U(x),\dim\ovl M_V(x)\in\{0,1\}$. We divide the argument into a case analysis.

\emph{Case 1.} $\dim\ovl M_U(x) = 0 = \dim\ovl M_V(x)$.

By Definition \ref{definition string}, we have that $\g^{-1}(x)\cap U = \emp = \g^{-1}(x)\cap V$. Therefore, by Lemma \ref{lemma I J}, $\g^{-1}(x)\cap I = \g^{-1}(x)\cap((U-2l\pi)\cup V) = ((\g^{-1}(x)\cap U)-2l\pi)\cup (\g^{-1}(x)\cap V) = \emp$ and $\g^{-1}(x)\cap J = \g^{-1}(x)\cap U\cap (V+2l\pi) = \emp$. Thus, $\dim\ovl M_I(x) = 0 = \dim\ovl M_J(x)$.

\emph{Case 2.} $\dim\ovl M_U(x) = 1$ and $\dim\ovl M_V(x) = 0$.

We have that $\g^{-1}(x)\cap U = \{s\}$ for some $s\in\R$, and $\g^{-1}(x)\cap V = \emp$. Therefore, $\g^{-1}(x)\cap I =  ((\g^{-1}(x)\cap U)-2l\pi)\cup (\g^{-1}(x)\cap V) = \{s-2l\pi\}$ and $\g^{-1}(x)\cap J = U\cap ((\g^{-1}(x)\cap V)+2l\pi) = \emp$. We obtain that $\dim\ovl M_I(x) = 1$ and $\dim \ovl M_J(x) = 0$.

\emph{Case 3.} $\dim\ovl M_U(x) = 0$ and $\dim\ovl M_V(x) = 1$.

The proof is similar to Case 2.

\emph{Case 4.} $\dim\ovl M_U(x) = 1 = \dim\ovl M_V(x)$.

In this case, $\g^{-1}(x)\cap U = \{s\}$ and $\g^{-1}(x)\cap V = \{t\}$ for some $s,t\in\R$. Thus, $\g^{-1}(x)\cap I =  ((\g^{-1}(x)\cap U)-2l\pi))\cup (\g^{-1}(x)\cap V) = \{s-2l\pi,t\}$ and $\g^{-1}(x)\cap J = (\g^{-1}(x)\cap U)\cap ((\g^{-1}(x)\cap V)+2l\pi) = \{s\}\cap\{t+2l\pi\}$. If $s = t+2l\pi$, then $\g^{-1}(x)\cap I = \{t\}$ and $\g^{-1}(x)\cap J = \{s\}$, i.e. $\dim\ovl M_I(x) = 1 = \dim\ovl M_J(x)$. If $s \neq t+2l\pi$, then $\dim\ovl M_I(x) = 2$ and $\dim\ovl M_J(x) = 0$. This concludes the argument of statement (1).

Now we prove statement (2). We check that the claim holds when $\dim\ovl M_U(x) = 1 = \dim\ovl M_V(x)$, the other cases are analogous to the above. Let $s,t\in\R$ be such that $\g^{-1}(x)\cap U = \{s\}$ and $\g^{-1}(x)\cap V = \{t\}$. If $s = t+2l\pi$, then $s\in (\g^{-1}(x)\cap U)\cap (\g^{-1}(x)\cap (V+2l\pi)) = \g^{-1}(x)\cap (U\cap(V+2l\pi)) = \g^{-1}(x)\cap J$. As a consequence, $J\neq\emp$, and this gives a contradiction. Thus, $s\neq t+2l\pi$, and we obtain that $\dim\ovl M_I(x) = 2$ similarly to the above.
\end{proof}

Now we prove that certain sequences of objects and morphisms of $\rep(\Zc_m,\ka)$ are exact.

\begin{proposition}\label{proposition short exact sequence}
Keeping Setup \ref{setup short exact sequences}, the following statements hold.
\begin{enumerate}
\item If $J \neq \emp$, the sequence $0 \lora \ovl M_U\ovs{\left(\begin{smallmatrix}f_1\\ f_2\end{smallmatrix}\right)}{\lora} \ovl M_I\oplus \ovl M_J\ovs{\left(\begin{smallmatrix} g_1 & g_2\end{smallmatrix}\right)}{\lora} \ovl M_V\lora 0$ with $f_1,f_2,g_1,g_2\neq 0$ such that $g_1f_1+g_2f_2 = 0$, is exact.
\item If $J = \emp$, the sequence $0 \lora \ovl M_U\ovs{f}{\lora} \ovl M_I\ovs{g}{\lora} \ovl M_V\lora 0$ with $f,g\neq 0$ such that $gf = 0$, is exact.
\end{enumerate}
\end{proposition}
\begin{proof}
We prove (1), statement (2) is analogous. It is enough to show that for each $x\in S^1$
\begin{equation}\label{equation short exact sequence}
\begin{tikzcd}[ampersand replacement = \&]
0\ar[r] \& \ovl M_U(x)\ar[r,"{\left(\begin{smallmatrix}f_1(x)\\ f_2(x)\end{smallmatrix}\right)}"] \& \ovl M_I(x)\oplus \ovl M_J(x)\ar[rr,"{\left(\begin{smallmatrix}g_1(x) & g_2(x)\end{smallmatrix}\right)}"] \& \&\ovl M_V(x)\ar[r] \& 0
\end{tikzcd}
\tag{E}
\end{equation}
is a short exact sequence of vector spaces. Let $x\in S^1$, since $\ovl M_U$ and $\ovl M_V$ are not projective, $\ovl M_U(x)$ and $\ovl M_V(x)$ are either equal to $0$ or isomorphic to $\K$. We assume that $\ovl M_U(x)\cong\K$ and $\ovl M_V(x)\cong\K$, for the remaining cases the proof is analogous. By Lemma \ref{lemma short exact sequence}, we have that either $\ovl M_I(x)\cong\K$ and $\ovl M_J(x)\cong\K$, or $\ovl M_I(x)\cong\K^2$ and $\ovl M_J(x) = 0$, we proceed with a case analysis. 

\emph{Case 1.} $\ovl M_I(x)\cong\K$ and $\ovl M_J(x)\cong \K$. 

By Lemma \ref{lemma I J} , $f_1$, $f_2$, $g_1$ and $g_2$ are either monomorphisms or epimorphisms. Thus, $f_1(x)$, $f_2(x)$, $g_1(x)$, and $g_2(x)$ are either injective or surjective linear maps, in particular they are non-zero maps between one-dimensional vector spaces. Then, (\ref{equation short exact sequence}) is isomorphic to the short exact sequence
\[
0 \lora \K\ovs{\left(\begin{smallmatrix}1\\ 1\end{smallmatrix}\right)}{\lora} \K\oplus\K\ovs{\left(\begin{smallmatrix} -1 & 1\end{smallmatrix}\right)}{\lora} \K\lora 0
\]
and is therefore is a short exact sequence. 

\emph{Case 2.} $\ovl M_I(x)\cong\K^2$ and $\ovl M_J(x) = 0$. 

We denote $f(x) = \left(\begin{smallmatrix}\l_1 & \l_2\end{smallmatrix}\right)^T$ and $g(x) = \left(\begin{smallmatrix}\mu_1 & \mu_2\end{smallmatrix}\right)$ for some $\l_1,\l_2,\mu_1,\mu_2\in\K$. We recall that, since $\ovl M_I(x)\cong\K^2$, $\ovl M_I$ is projective and $\ovl M_I(\om_x)\colon \ovl M_I(x)\to \ovl M_I(x)$ is of the form $\ovl M_I(x) = \left(\begin{smallmatrix}0 & 0\\ 1 & 0\end{smallmatrix}\right)$, see Definition \ref{definition string}. Since $\ovl M_I(\om_x)f(x) = f(x)\ovl M_U(\om_x) = 0$, we obtain that $f(x) = \left(\begin{smallmatrix} 0 & \l_2\end{smallmatrix}\right)^T$. Moreover, since $f(x)\neq 0$ as in Case 1, we have that $\l_2\neq 0$. Similarly, we obtain that $g(x) = \left(\begin{smallmatrix}\mu_1 & 0\end{smallmatrix}\right)$ with $\mu_1\neq 0$. Thus, (\ref{equation short exact sequence}) is isomorphic to the short exact sequence
\[
0 \lora \K\ovs{\left(\begin{smallmatrix}0\\ 1\end{smallmatrix}\right)}{\lora} \K\oplus\K\ovs{\left(\begin{smallmatrix} 1 & 0\end{smallmatrix}\right)}{\lora} \K\lora 0.
\]
We obtain that (\ref{equation short exact sequence}) is a short exact sequence for each $x\in S^1$. 
\end{proof}

\section{Simple objects and uniseriality}\label{section simple objects rep}

With the following proposition, which is analogous to \cite[Remark 6.1(1)]{RZ}, we describe the simple objects of $\rep(\Zc_m,\ka)$. For each $z\in\Zc_m$ we denote $S_z = \ovl M_{(z,z^+]}$.

\begin{proposition}\label{proposition simple objects}
The simple objects of $\rep(\Zc_m,\ka)$ are exactly those of the form $S_z$ for $z\in\Zc_m$.
\end{proposition}
\begin{proof}
Let $z\in\Zc_m$, we show that $S_z$ is simple. Let $U = (u_1,u_2+2h\pi]\in\I$ be such that there exists a monomorphism $0\to \ovl M_U\to S_z$. By Lemma \ref{lemma mono and epi} we have that $((z,z^+]+2l\pi)\cap_L U\neq\emp$ and $z^++2l\pi = u_2+2h\pi$ for a unique $l\in\{0,-1\}$. Since $(z,z^+]\subseteq(0,2\pi)$, then $l = 0$ and as a consequence $h = 0$ and $u_2 = z^+$. Since $(z,z^+]\cap_L U\neq\emp$, we obtain that $U = (z,z^+]$. This proves that $S_z$ is simple.

Now let $V = (v_1,v_2+2k\pi]\in\I$ and assume that $\ovl M_V$ is simple. Note that, by Lemma \ref{lemma mono and epi}, there exists a monomorphism $0\to S_{v_2^-}\to \ovl M_V$. Thus, $S_{v_2^-}\cong\ovl M_V$, i.e. $V = (v_2^-,v_2]$. This concludes the argument.
\end{proof}

From the lemma below it follows that for each $U = (u_1,u_2+2h\pi]\in\I$, the object $\ovl M_U$ admits the following, possibly countable infinite, composition series
\[
0\subseteq\ovl M_{(u_2^-,u_2]}\subseteq\ovl M_{(u_2^{--},u_2]} \subseteq \cdots \subseteq\ovl M_{(u_1^{++},u_2+2h\pi]}\subseteq \ovl M_{(u_1^+,u_2+2h\pi]}\subseteq\ovl M_{(u_1,u_2+2h\pi]} = \ovl M_U.
\]

\begin{lemma}
For each $M\in\ind\rep(\Zc_m,\ka)$ there exists a monomorphism $f\colon L\to M$ such that $\Coker f$ is simple. Moreover, $f$ is unique up to isomorphism.
\end{lemma}
\begin{proof}
Assume that $M$ is simple, then $0\to M$ is the only monomorphism whose cokernel is simple. Now, assume that $M$ is not simple, and let $U = (u_1,u_2+2h\pi]\in\I$ be such that $M\cong\ovl M_U$. Let $V = (u_1^+,u_2+2h\pi]$, note that $V\in\I$ and, by Lemma \ref{lemma mono and epi}, there exists a monomorphism $f\colon \ovl M_V\to \ovl M_U$. Moreover, by Proposition \ref{proposition short exact sequence} , $0\lora \ovl M_V\ovs{f}{\lora}\ovl M_U\lora S_{u_1}\lora 0$ is a short exact sequence, i.e. $\Coker f$ is simple. 

Now we prove the uniqueness. Assume that there exists $V = (v_1,v_2+2k\pi]\in\I$ and a monomorphism $f\colon \ovl M_V\to \ovl M_U$ such that $\Coker f$ is simple. Assume that $h = 0$, the other case is similar. By Lemma \ref{lemma mono and epi}, $k = 0$ and $v_2 = u_2$. Moreover, by Proposition \ref{proposition short exact sequence}, it is straightforward to check that $\Coker f \cong \ovl M_{(u_1,v_1]}$. Thus, $v_1 = u_1^+$ and $V = (v_1,v_2] = (u_1^+,u_2]$. This concludes the proof.
\end{proof}

We obtain that the category $\rep(\Zc_m,\ka)$ is uniserial in the following sense.

\begin{theorem}\label{theorem uniserial}
Each indecomposable object of $\rep(\Zc_m,\ka)$ has a unique, possibly countable infinite, composition series.
\end{theorem}

\section{Irreducible morphisms and almost split sequences}\label{section irreducible morphisms and almost split sequences in rep}

We describe the irreducible morphisms and the almost split sequences in $\rep(\Zc_m,\ka)$. Given $U\in\I$, we refer to Definition \ref{definition intervals} for the intervals $U_1$, $U_2$, and $U^-$. 

\begin{proposition}\label{proposition irreducible morphism in rep}
Let $U\in\I$. If $\ovl M_U$ is not projective, consider a non-zero morphism $f_1\colon \ovl M_U\to \ovl M_{U_1}$. If $\ovl M_U$ is not simple, consider a non-zero morphism $f_2\colon \ovl M_U\to \ovl M_{U_2}$. The following statements hold.
\begin{enumerate}
\item Let $V\in\I$, $V\neq U$. Then any morphism $\ovl M_U\to \ovl M_V$ factors through $f_1$ or $f_2$.
\item If $\ovl M_U$ is not simple, then $f_2$ is irreducible.
\item If $\ovl M_U$ is not projective, then $f_1$ is irreducible.
\item Let $V\in\I$ be such that that there exists an irreducible morphism $\ovl M_U\to \ovl M_V$. Then $\ovl M_V\cong \ovl M_{U_1}$ or $\ovl M_V\cong\ovl M_{U_2}$.
\end{enumerate}
\end{proposition}
\begin{proof}
We prove statement (1). Assume that $\ovl M_U$ is not projective and not simple, for the other cases the proof is analogous. By Observation \ref{observation intervals}, $U_1, U_2\in\I$ and it is straightforward to check that $U_1\cap_L U\neq\emp$ and $U_2\cap_L U\neq\emp$. Now, we write $U = (u_1,u_2+2h\pi]$, $V = (v_1,v_2+2k\pi]$, and we consider $f\colon\ovl M_U\to \ovl M_V$. If $f = 0$ then $f$ factors through $f_1$ and $f_2$. We assume that $f \neq 0$. If $v_1 = u_1$, then $V\cap_L U\neq\emp$, $V\cap_L U_2\neq\emp$ and, by Proposition \ref{proposition factorisation properties in rep}, $h$ factors through $f_2$. If $v_1\neq u_1$, then we have the following possibilities: either $V\cap_L U\neq\emp$, and then $V\cap_L U_1\neq\emp$, or $(V-2\pi)\cap_L U\neq\emp$, and then $(V-2\pi)\cap_L U_1\neq\emp$. In both cases we have that $f$ factors through $f_1$. We conclude that $f$ factors through $f_1$ or $f_2$.

Now we prove statement (2), the proof of statement (3) is analogous. Since $\ovl M_U$ is not simple we have that $U_2\in\I$. Assume that $f_2\colon \ovl M_U\to \ovl M_{U_2}$ factors as $f_2 = \b\a$ for some $\a\colon \ovl M_U\to M$, $\b\colon M\to \ovl M_{U_2}$, and $M\in\rep(\Zc_m,\ka)$. We show that $\a$ is a split monomorphism or $\b$ is a split epimorphism. Since $\Hom_{S^1}(\ovl M_U,\ovl M_{U_2})\cong\K$, without loss of generality we assume that $M$ is indecomposable, i.e. $M\cong \ovl M_V$ for some $V\in\I$. Since $\b\a = f_2$ and $U_2\cap_L U\neq\emp$, we have that $V\cap_L U\neq \emp$ and $U_2\cap_L V\neq\emp$. Indeed, for the remaining cases, by Proposition \ref{proposition zero compositions in rep}, we have that $\b\a = 0$, giving a contradiction. Thus, $V = U$ or $V = U_2$, i.e. $\a$ is a split monomorphism or $\b$ is a split epimorphism. We conclude that $f_2$ is irreducible.

We prove statement (4). Consider an irreducible morphism $f\colon \ovl M_U\to \ovl M_V$, by statement (1) $f$ factors through $f_1$ or $f_2$. Assume that $f$ factors through $f_1$, the other case is analogous. Thus, there exists $g\colon \ovl M_{U_1}\to \ovl M_V$ such that $f =gf_1$. Since $f$ is irreducible and $f_1$ is not a split monomorphism, $g$ is a split epimorphism and therefore $\ovl M_{U_1}\cong \ovl M_V$.
\end{proof}

\begin{proposition}\label{proposition almost split sequences in rep}
Let $U\in\I$ be such that $\ovl M_U$ is not projective. The following statements hold.
\begin{enumerate}
\item If $\ovl M_U$ is not simple, the short exact sequence $0\lora \ovl M_U\lora \ovl M_{U_1}\oplus\ovl M_{U_2}\lora \ovl M_{U^-}\lora 0$ is almost split.	
\item If $\ovl M_U$ is simple, the short exact sequence $0\lora \ovl M_U\lora \ovl M_{U_1}\lora \ovl M_{U^-}\lora 0$ is almost split.
\end{enumerate}
Thus, we have that $\t\ovl M_{U^-} \cong \ovl M_U$.
\end{proposition}
\begin{proof}
We prove statement (1), statement (2) is analogous. We denote $\left(\begin{smallmatrix}f_1 & f_2\end{smallmatrix}\right)^T\colon \ovl M_U\to \ovl M_{U_1}\oplus \ovl M_{U_2}$. The objects $\ovl M_U$ and $\ovl M_{U^-}$ are indecomposable, thus it remains to check that $\left(\begin{smallmatrix}f_1 & f_2\end{smallmatrix}\right)^T$ is left almost split. Consider a morphism $\a\colon \ovl M_U\to M$ with $M\in\rep(\Zc_m,\ka)$, and assume that $\a$ is not a split monomorphism, we show that $\a$ factors through $\left(\begin{smallmatrix}f_1 & f_2\end{smallmatrix}\right)^T$. Without loss of generality we can assume that $M$ is indecomposable. By Proposition \ref{proposition irreducible morphism in rep}, $\a$ factors through $f_1$ or $f_2$. Assume that $\a = \a_1f_1$ for some $\a_1\colon \ovl M_{U_1}\to M$, the other case is similar. Consider the morphism $\left(\begin{smallmatrix}\a_1 & 0\end{smallmatrix}\right)\colon \ovl M_{U_1}\oplus \ovl M_{U_2}\to M$, we obtain that $\left(\begin{smallmatrix}\a_1 & 0\end{smallmatrix}\right)\left(\begin{smallmatrix}f_1 & f_2\end{smallmatrix}\right)^T = \a$. Thus, we obtain that $\left(\begin{smallmatrix}f_1 & f_2\end{smallmatrix}\right)^T$ is left almost split. 
\end{proof}

\section{A negative CY triangulated category}\label{section a negative calabi-yau triangulated category}

We define the category $\C_{-1,m}$ and we describe its geometric model and AR quiver. We prove that $\C_{-1,m}$ is $(-1)$-CY, and we observe that for the case $m = 1$ it is triangle equivalent to the $(-1)$-CY version of the Holm--J{\o}rgensen category, which first appeared in \cite{HJY}. 

\begin{definition}
The category $\C_{-1,m}$ is the stable category of $\rep(\Zc_m,\ka)$.
\end{definition}

Since $\rep(\Zc_m,\ka)$ is a Frobenius category, then $\C_{-1,m}$ is triangulated. Moreover, $\C_{-1,m}$ is a $\K$-linear, $\Hom$-finite, Krull--Schmidt category because $\rep(\Zc_m,\ka)$ is. We refer to Section \ref{section frobenius categories} for some background on stable Frobenius categories.

\subsection{Geometric model}\label{section the geometric model -1 cy}

The indecomposable objects of $\C_{-1,m}$ are the non-projective-injective indecomposable objects of $\rep(\Zc_m,\ka)$. More precisely, they are of the form $\ovl M_U$ with $U = (u_1,u_2+2h\pi]\in\I$ such that $u_1^+\leq u_2+2h\pi\leq u_1+2\pi$. 

The $\Hom$-spaces of $\C_{-1,m}$ are at most one-dimensional because the $\Hom$-spaces of $\rep(\Zc_m,\ka)$ between non-projective-injective objects are, see Proposition \ref{proposition hom spaces}. We recall that the morphisms of $\C_{-1,m}$ are the equivalence classes $[f]$ of morphisms $f$ of $\rep(\Zc_m,\ka)$, we refer to Section \ref{section frobenius categories} for more details. In order to simplify the notation, we omit the square brackets in this section. It will be clear from the context when we refer to morphisms of $\rep(\Zc_m,\ka)$ or $\C_{-1,m}$.

The action of the shift functor $\S\colon\C_{-1,m}\to \C_{-1,m}$ on objects is as follows. Let $U\in\I$ be such that $\ovl M_U\in\rep(\Zc_m,\ka)$ is an indecomposable non-projective-injective object. By Proposition \ref{proposition short exact sequence}, the following is a short exact sequence in $\rep(\Zc_m,\ka)$, and therefore $\S\ovl M_U = \ovl M_{\S U}$.
\[
0\lora \ovl M_U\lora I_{u_2^-}\lora \ovl M_{\S U}\lora 0
\]

We describe the geometric model of $\C_{-1,m}$ in terms of the $\infty$-gon $\Zc_m$ introduced in Section \ref{section the ifinity-gon}. Given $a_1,a_2\in\Zc_m$, when we write $a_1-a_2$ we mean the difference of $a_1$ and $a_2$ regarded as integers, forgetting about what copy of $\Z$ they belong to: if $a_1\in\Z^{(p)}$ and $a_2\in\Z^{(q)}$ for $p,q\in[m]$, then we write $a_1 = (n_1,p)$ and $a_2 = (n_2,q)$ for $n_1,n_2\in\Z$, and we set $a_1-a_2 = n_1-n_2$.

The following is a ``negative version" of \cite[Definition 2.3]{HJ2}. 

\begin{definition}
Let $a_1,a_2\in\Zc_m$. We say that $(a_1,a_2)$ is an \emph{admissible arc} if $a_1\geq a_2+1$ and $a_1-a_2\equiv 1\mod 2$.
\end{definition}

We prove that the admissible arcs of $\Zc_m$ are in bijection with certain intervals of $\I$. 

\begin{proposition}\label{proposition bijection intervals and -1-admissible arcs}
The following is a bijection.
\begin{align*}
\phi\colon \left\{\parbox{6cm}{\centering $U = (u_1,u_2+2h\pi]\in\I$ such that $u_1^+\leq u_2+2h\pi\leq u_1+2\pi$} \right\}& \longrightarrow \left\{\parbox{3.75cm}{\centering Admissible arcs of $\Zc_m$} \right\}\\
(u_1,u_2+2h\pi] & \longmapsto 
\begin{cases}
(2u_2-1,2u_1) & \text{if $h = 0$,}\\
(2u_1,2u_2-1) & \text{if $h = 1$.}
\end{cases}
\end{align*}
\end{proposition}
\begin{proof}
It is straightforward to check that $\phi$ is well defined and injective. We show that $\phi$ is surjective. Let $(x_1,x_2)$ be an admissible arc and let $U = \left(\frac{x_1}{2}, \frac{x_2+1}{2}+2\pi\right]$ if $x_1$ is even, and $U = \left(\frac{x_2}{2}, \frac{x_1+1}{2}\right]$ if $x_1$ is odd. Then $U$ belongs to the domain of $\phi$ and $\phi(U) = (x_1,x_2)$.
\end{proof}

It is straightforward to check that the functor $\S$ acts as a clockwise rotation of admissible arcs, i.e. for $a = (a_1,a_2)\in\ind\C_{-1,m}$ we have that $\S a = (a_1-1,a_2-1)$. For the rest of this section we identify the indecomposable objects of $\C_{-1,m}$ with the admissible arcs.

\begin{figure}[h]
\centering
\includegraphics[height = 5cm]{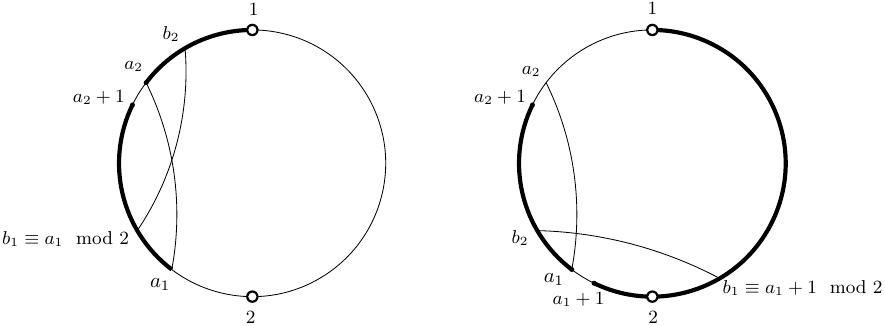}
\caption{The sets $H^+(a)$, on the left, and $H^-(\S^{-1}a)$, on the right, for $a\in\ind\C_{-1,2}$.}
\label{figure morphisms geometric model}
\end{figure}

For an object $a = (a_1,a_2)\in\ind\C_{-1,m}$, we re-write the $\Hom$-hammocks of Definition \ref{definition hom hammocks} as
\begin{align*}
H^+(a) & = \{(b_1,b_2)\in\ind\C_{-1,m}\mid a_2+1\leq b_1\leq a_1,   b_2\leq a_2 \text{, and } b_1\equiv a_1\mod 2\},\text{ and}\\
H^-(a) & = \{(b_1,b_2)\in\ind\C_{-1,m}\mid b_1\geq a_1,  a_2\leq b_2\leq a_1-1 \text{, and } b_1\equiv a_1\mod 2\}.
\end{align*}
By Proposition \ref{proposition hom hammocks} it follows that
\[
\Hom_{\C_{-1,m}}(a,b) \cong 
\begin{cases}
\K & \text{if $b\in H^+(a)\cup H^-(\S^{-1}a)$,}\\  
0 & \text{otherwise.}
\end{cases}
\]
Figures \ref{figure morphisms geometric model} and \ref{figure hom-hammocks} illustrate, respectively, the sets $H^+(a)$ and $H^-(\S^{-1}a)$ in the geometric model and in the AR quiver of $\C_{-1,m}$. 

\subsection{AR quiver}\label{section AR quiver}

For each $p,q\in[m]$ with $p\geq q$ and $i\in\{0,1\}$, we define the set of admissible arcs
\[
\Z^{(p,q,i)} = \left\{a = (a_1,a_2)\text{ admissible arc}\mid a_1\in\Z^{(p)}, a_2\in\Z^{(q)},\text{ and } a_1\equiv i\mod 2\right\}.
\]
We can arrange the isoclasses of indecomposable objects of $\C_{-1,m}$ into a coordinate system having
\begin{itemize}
\item $2m$ components of type $\Z A_{\infty}$, each corresponding to the sets of arcs $\Z^{(p,p,i)}$ for $p\in[m]$ and $i\in\{0,1\}$,
\item $2\binom{m}{2}$ components of type $\Z A_{\infty}^{\infty}$, each corresponding to the sets of arcs $\Z^{(p,q,i)}$ for $p> q$ and $i\in\{0,1\}$.
\end{itemize}

\begin{figure}[h]
\centering
\includegraphics[height = 4cm]{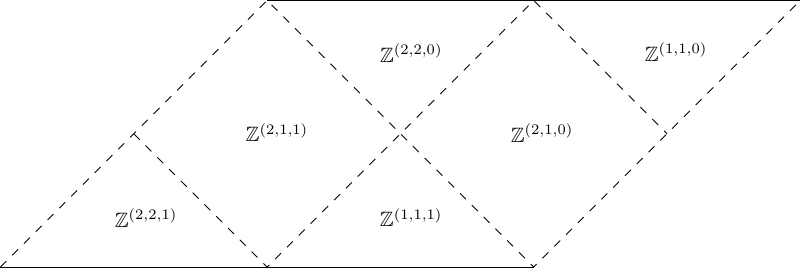}
\caption{The coordinate system of $\C_{-1,2}$.}
\label{figure ar quiver}
\end{figure}

Figure \ref{figure ar quiver} provides an illustration of the coordinate system, with Proposition \ref{proposition irreducible morphisms} we will prove that this yields the AR quiver of $\C_{-1,m}$. 

\begin{figure}[h]
\centering
\includegraphics[height = 5cm]{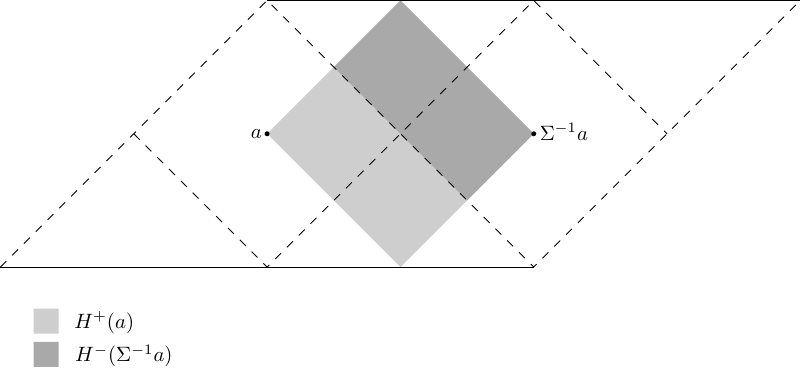}
\caption{The $\Hom$-hammocks $H^+(a)$ and $H^-(\S^{-1}a)$ for $a\in\ind\C_{-1,2}$.}
\label{figure hom-hammocks}
\end{figure}

Figure \ref{figure hom-hammocks} provides an illustration of the $\Hom$-hammocks in the coordinate system.

\subsection{CY property}\label{section -1 cy}

We prove that the category $\C_{-1,m}$ is $(-1)$-CY, i.e. $\S^{-1}$ is a Serre functor. By \cite[Proposition I.1.4]{RV} it is enough to prove that for each $a,b\in\C_{-1,m}$ there exists a \emph{non-degenerate pairing} $\Phi_{a,b}\colon \Hom(a,b)\times\Hom(b,\S^{-1}a)\to \K$, i.e. a bilinear map such that if $\Phi_{a,b}(f,g) = 0$ for each $f\in \Hom(a,b)$ then $g = 0$, and if $\Phi_{a,b}(f,g) = 0$ for each $g\in \Hom(\S^{-1}b,a)$ then $f = 0$.

We now introduce some notation. Consider $a\in\C_{-1,m}$, its decomposition into indecomposable direct summands $a_1\oplus\cdots\oplus a_n$, and a morphism $f\colon a\to \S^{-1}a$. We can write $f$ as the matrix $f = (f_{i,j})_{i,j}$, where $f_{i,j}\colon a_j\to \S^{-1}a_i$. We indicate the \emph{trace} of $f$ by $\Tr(f) = f_{1,1}+\cdots+f_{n,n}$.

\begin{theorem}\label{theorem non degenerate bilinear form}
Let $a,b\in\C_{-1,m}$. The following is a non-degenerate pairing.
\begin{align*}
\Phi_{a,b}\colon \Hom_{\C_{-1,m}}(a,b)\times \Hom_{\C_{-1,m}}(b,\S^{-1}a)& \lora \K\\
(f, g) & \longmapsto \Tr(gf)
\end{align*}
\end{theorem}

\begin{corollary}\label{corollary -1-cy}
The category $\C_{-1,m}$ is $(-1)$-CY, i.e. $\S^{-1}$ is a Serre functor.
\end{corollary}

We will prove Theorem \ref{theorem non degenerate bilinear form} after we discuss the factorisation properties of the morphisms of $\C_{-1,m}$. We refer to Proposition \ref{proposition factorisation properties in rep} for the factorisation properties in $\rep(\Zc_m,\ka)$.

\begin{proposition}\label{proposition factorisation properties}
Let $a,b,c\in\ind\C_{-1,m}$, $f\colon a\to b$, and $g\colon b\to c$ be non-zero morphisms. Assume that one of the following conditions holds.
\begin{enumerate}
\item $b\in H^+(a)$ and $c\in H^+(a)\cap H^+(b)$.
\item $b\in H^+(a)$ and $c\in H^-(\S^{-1}a)\cap H^-(\S^{-1}b)$.
\item $b\in H^-(\S^{-1}a)$ and $c\in H^-(\S^{-1}a)\cap H^+(b)$.
\end{enumerate}
Then $gf\neq 0$.
\end{proposition}
\begin{proof}
We recall that $\phi$ denotes the bijection of 	Proposition \ref{proposition bijection intervals and -1-admissible arcs}. Let $U = (u_1,u_2+2h\pi], V = (v_1,v_2+2k\pi], W = (w_1,w_2+2l\pi]\in\I$ be intervals such that $\phi(U) = a$, $\phi(V) = b$, and $\phi(W) = c$. Assume that (1) holds, i.e. $V\in H^+(U)$ and $W\in H^+(U)\cap H^+(V)$, then $V\cap_L U\neq \emp$, $W\cap_L V\neq\emp$, and $W\cap_L U\neq\emp$, see Observation \ref{observation h+ h- p}. Thus, by Proposition \ref{proposition factorisation properties in rep}, $gf\neq 0$ in $\rep(\Zc_m,\ka)$. Since $W\notin P(U')$, $gf\neq 0$ in $\C_{-1,m}$. 

Now assume that (2) holds, i.e. $V\in H^+(U)$ and $W\in H^-(\S^{-1}U)\cap H^-(\S^{-1}V)$. By Observation \ref{observation h+ h- p}, if $h = 0$ then $k = h = 0$, $l = 1-k = 1$, $V\cap_L U\neq\emp$, $(W-2\pi)\cap_L V\neq\emp$, and $(W-2\pi)\cap_L U\neq\emp$. Thus, $gf\neq 0$. If $h = 1$, then $k = h = 1$, $l = 1-k = 0$, $V\cap_L U\neq\emp$, $W\cap_L V\neq\emp$, and $W\cap_L U\neq\emp$. Similarly as above, we obtain that $gf\neq 0$ in $\C_{-1,m}$.

Finally, assume that (3) holds. If $h = 0$, then $k = 1-h = 1$, $l = k = 1$, $(V-2\pi)\cap_L U\neq\emp$, $W\cap_L V\neq\emp$, and $(W-2\pi)\cap_L U\neq\emp$. If $h = 1$, then $k = 1-h = 0$, $l = k = 0$, $V\cap_L U\neq\emp$, $W\cap_L V\neq\emp$, and $W\cap_L U\neq\emp$. In both cases we have that $gf\neq 0$ in $\C_{-1,m}$.
\end{proof} 

The following lemma is useful for proving Theorem \ref{theorem non degenerate bilinear form}. 

\begin{lemma}\label{lemma -1 cy}
Let $a,b\in\ind\C_{-1,m}$ and $f\colon a\to b$ be non-zero. Then there exists $g\colon b\to \S^{-1}a$ such that $gf\neq 0$.
\end{lemma}
\begin{proof}
Since $f\neq 0$, $b\in H^+(a)\cup H^-(\S^{-1}a)$. If $b\in H^+(a)$ then $a\in H^-(b)$, i.e. $\S^{-1}a\in H^-(\S^{-1}b)$, and if $b\in H^-(\S^{-1}a)$ then $\S^{-1}a\in H^+(b)$. Thus, there exists a non-zero morphism $g\colon b\to \S^{-1}a$. We have the following possibilities: $b\in H^+(a)$ and $\S^{-1}a\in H^-(\S^{-1}a)\cap H^-(\S^{-1}b)$, or $b\in H^-(\S^{-1}a)$ and $\S^{-1}a\in H^-(\S^{-1}a)\cap H^+(b)$. In both cases, by Proposition \ref{proposition factorisation properties}, we obtain that $gf\neq 0$.
\end{proof}

\begin{proof}[{Proof of Theorem \ref{theorem non degenerate bilinear form}}]
It is straightforward to check that $\Phi_{a,b}$ is a bilinear map, we prove that it is non-degenerate. Let $f\colon a\to b$ be non-zero, we show that there exists $g\colon b\to \S^{-1}a$ such that $\Phi_{a,b}(f,g) = \Tr(gf)\neq 0$. We denote $a = \bigoplus_{t = 1}^n a_t$ and $b =  \bigoplus_{s = 1}^k b_s$ with $a_1,\dots,a_n,b_1,\dots,b_k\in\ind\C_{-1,m}$, and $f = (f_{s,t})_{s,t}$ with $f_{s,t}\colon a_t\to b_s$. Since $f\neq 0$, there exist $1\leq i\leq n$ and $1\leq j\leq k$ such that $f_{i,j}\colon a_j\to b_i$ is non-zero. By Lemma \ref{lemma -1 cy}, there exists $\g\colon b_i\to \S^{-1}a_j$ such that $\g f_{i,j}\neq 0$. Now, we define $g = (g_{s,t})_{s,t}\colon b\to \S^{-1}a$ as
\[
g_{s,t} = 
\begin{cases}
\g & \text{if $s = j$ and $t = i$,}\\
0 & \text{otherwise.}
\end{cases}
\]
Then, $(gf)_{t,t} = \g f_{i,j}\neq 0$ if $t = j$, and $(gf)_{t,t} = 0$ otherwise. As a consequence, $\Tr(gf)\neq 0$. Similarly, we can prove that given a non-zero morphism $g\colon b\to \S^{-1}a$, there exists $f\colon a\to b$ such that $\Phi_{a,b}(gf) = \Tr(gf)\neq 0$. We conclude that $\Phi_{a,b}$ is non-degenerate.
\end{proof}

\subsection{Irreducible morphisms and almost split sequences}

We describe the irreducible morphisms and the almost split sequences in $\C_{-1,m}$. The following result follows directly from Propositions \ref{proposition irreducible and left right almost split} and \ref{proposition irreducible morphism in rep}.

\begin{proposition}\label{proposition irreducible morphisms}
Let $a = (a_1,a_2),b = (b_1,b_2)\in\ind\C_{-1,m}$. If $b = (a_1,a_2-2)$ or $b = (a_1-2,a_2)$, then any non-zero morphism $a\to b$ of $\C_{-1,m}$ is irreducible. Moreover, there are no other irreducible morphisms in $\C_{-1,m}$ between indecomposable objects.
\end{proposition}

The following result follows directly from Corollary \ref{corollary almost split triangles after stabilising} and Proposition \ref{proposition almost split sequences in rep}.

\begin{proposition}\label{proposition almost split triangles}
Let $a = (a_1,a_2)\in\ind\C_{-1,m}$. The following statements hold. 
\begin{enumerate}
\item If $a_1 = a_2+1$, then $(a_1+2,a_2+2)\lora(a_1+2,a_2)\lora (a_1,a_2)\lora \S(a_1+2,a_2+2)$ is an almost split triangle.
\item If $a_1\neq a_2+1$, then $(a_1+2,a_2+2)\lora(a_1,a_2+2)\oplus (a_1+2,a_2)\lora (a_1,a_2)\lora \S(a_1+2,a_2+2)$ is an almost split triangle.
\end{enumerate}
Thus, we have that $\t a = (a_1+2,a_2+2)$.
\end{proposition}

\begin{figure}[h]
\centering
\includegraphics[height = 5cm]{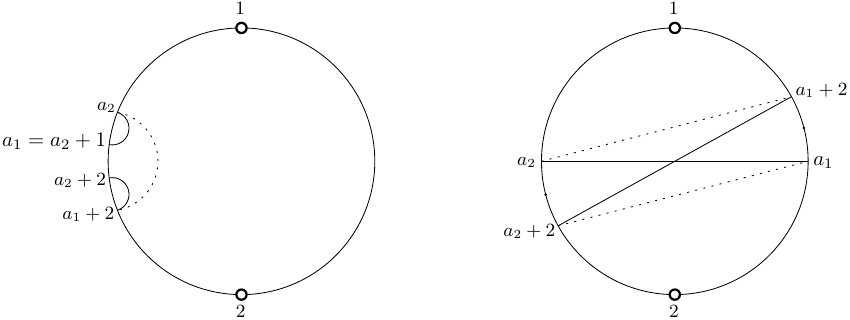}
\caption{The almost split triangles in $\C_{-1,2}$. The dotted arcs represent the middle terms of the triangles.}
\label{figure almost split triangles}
\end{figure}

Figure \ref{figure almost split triangles} provides an illustration of the almost split triangles in $\C_{-1,m}$.

\subsection{The case $m = 1$}

In Sections \ref{section the geometric model -1 cy} and \ref{section AR quiver} we described the geometric model and AR quiver of $\C_{-1,m}$. We observe that, for the case $m = 1$, these coincide with the geometric model and AR quiver of $\T_{-1}$, which denotes the $(-1)$-CY version of the Holm--J{\o}rgensen category introduced in \cite{HJ1, J}. We refer to \cite{C1, HJY} for the definition, geometric model, and AR quiver of $\T_{-1}$. We want to prove the following result.

\begin{theorem}\label{theorem triangle equivalence HJ category}
The categories $\C_{-1,1}$ and $\T_{-1}$ are triangle equivalent.
\end{theorem}

We start by recalling the following definition (see for instance \cite[Section 2]{CP1}). Given a $\K$-linear, $\Hom$--finite, Krull--Schmidt triangulated category $\T$, an object $t\in\T$ is called \emph{$(-1)$-spherical} if the following conditions hold.

\begin{itemize}
\item $\Hom(t,\S^n t)\cong \K$ if $n = 0$ or $n = -1$, and $\Hom(t,\S^n t) = 0$ otherwise.
\item $\Hom(t,u)\cong D\Hom(\S^{-1} u,t)$ for each $u\in\T$ and this isomorphism is natural in $u$. Here $D$ denote the usual vector-space duality $D = \Hom_{\K}(-,\K)$. 
\end{itemize}
The category $\T_{-1}$ is algebraic triangulated and is classically generated by a $(-1)$-spherical object. By \cite[Theorem 2.1]{KYZ}, $\T_{-1}$ is unique with these properties. 

We recall that there exists a family of $w$-CY analogues of the Holm-J{\o}rgensen category for $w\in\Z$, each classically generated by a $w$-spherical object. We refer to \cite{CP1, HJ2, HJY} for more details.

Theorem \ref{theorem triangle equivalence HJ category} follows directly from the lemma below.

\begin{lemma}
The category $\C_{-1,1}$ is classically generated by a $(-1)$-spherical object.
\end{lemma}
\begin{proof}
Let $a = (1,0)\in\ind\C_{-1,1}$. By the description of the $\Hom$-spaces of $\C_{-1,m}$ in Section \ref{section the geometric model -1 cy}, it is straightforward to check that $\Hom(a,\S^n a)\cong \K$ if $n = 0$ or $n = -1$, and $\Hom(a,\S^n a) = 0$ otherwise. Moreover, from Corollary \ref{corollary -1-cy}, it follows that $a$ is a $(-1)$-spherical object. We denote by $\thick(a)$ the smallest thick subcategory of $\C_{-1,m}$ containing the object $a$. We prove that $\thick(a) = \C_{-1,m}$. 

We observe that each indecomposable object of $\C_{-1,m}$ is of the form $(n,n-k)$, for some $n\in\Z$ and $k\geq 1$ odd. We fix $n\in\Z$, and we prove that $(n,n-k)\in\thick(a)$ by induction on $k$. For the case $k = 1$ we have that $\S^{-n+1}a = (n,n-1)\in\thick(a)$. If $k = 3$, then, by Proposition \ref{proposition almost split triangles}, $(n,n-k) = (n,n-3)$ is the middle term of the almost split triangle $(n,n-1)\lora (n,n-3) \lora (n-2,n-3)\lora \S(n,n-1)$. Since $(n,n-1)\in\thick(a)$ and $\S^2(n,n-1) = (n-2,n-3)\in\thick(a)$, then $(n,n-3)\in\thick(a)$. If $k\geq 5$, then $(n,n-k)$ is a direct summand of the middle term of the almost split triangle $(n,n-k+2) \lora (n-2,n-k+2)\oplus (n,n-k)\lora (n-2,n-k)\lora \S(n,n-k+2)$. By the induction hypothesis, $(n,n-k+2)\in\thick(a)$ and $\S^2(n,n-k+2) = (n-2,n-k)\in\thick(a)$, and as a consequence $(n-2,n-k+2)\in\thick(a)$. We conclude that $\thick(a) = \C_{-1,1}$. 
\end{proof}

\end{document}